\documentclass[11pt]{amsart}

\usepackage{stmaryrd}
\usepackage[linesnumbered,ruled,vlined]{algorithm2e}
\usepackage{subcaption}

\newtheorem{remark}{Remark}
\newtheorem{definition}{Definition}

\usepackage{hyperref}

\newcommand\EE{\mathbb E}

\newcommand\PP{\mathbb P}

\newcommand\RR{\mathbb R}

\newcommand\bJ{{\mathbf J}}
\newcommand\bN{\mathbf N}

\newcommand\bW{\mathbf W}

\newcommand{\cF}{{\mathcal F}}

\newcommand{\cL}{{\mathcal L}}

\newcommand{\cP}{{\mathcal P}}

\newcommand{\ctrl}{\alpha}
\newcommand{\ctrldim}{k}
\newcommand{\dom}{\mathcal{Q}}
\newcommand{\domT}{\mathcal{Q}_{T}}

\newcommand{\grad}{\nabla}

\DeclareMathOperator*{\bigtimes}{\vartimes}

\usepackage{empheq} 
\renewcommand{\div}{\mathrm{div}}

\begin{document}

\title[Deep Learning for MFG and MFC]{Deep Learning for Mean Field Games and Mean Field Control 
with Applications to Finance}
\author{Ren\'e Carmona \& Mathieu Lauri\`ere}
\address{Department of Operations Research and Financial Engineering, Princeton University, Email: rcarmona@princeton.edu, lauriere@princeton.edu. Work supported by NSF grant DMS-1716673 and ARO grant  W911NF-17-1-0578.}

\maketitle

\begin{abstract}
Financial markets and more generally macro-economic models involve a large number of individuals interacting through variables such as prices resulting from the aggregate behavior of all the agents. Mean field games have been introduced to study Nash equilibria for such problems in the limit when the number of players is infinite. The theory has been extensively developed in the past decade, using both analytical and probabilistic tools, and a wide range of applications have been discovered, from economics to crowd motion. More recently the interaction with machine learning has attracted a growing interest. This aspect is particularly relevant to solve very large games with complex structures, in high dimension or with common sources of randomness. In this chapter, we review the literature on the interplay between mean field games and deep learning, with a focus on three families of methods. A special emphasis is given to financial applications.

\end{abstract}

\tableofcontents

\section{Introduction}
Most applications in financial engineering rely on numerical implementations which suffer from the curse of dimensionality.
Recent developments in machine learning and the availability of powerful and readily available public domain packages have triggered a renewal of interest in financial numerics: deep learning technology has pushed the limits of computational finance, and one can now tackle problems which seemed out of reach a few years ago.

On a different part of the spectrum of scientific research, several sub-fields of economics experienced significant transitions: the emergence of continuous time stochastic models in macro-economics, and the development of general equilibrium theory in finance created a commonality between the two fields. This convergence provided a fertile ground for mean field game and mean field control theories, which appeared naturally as the tools of choice for theoreticians and applied mathematicians and economists.

Fashions come and go, and there is no point in trying to catch a train which already left the station. However, machine learning and mean field theories are here to stay. Gaining a deep understanding of the inner workings of both  paradigms  is the first step in the recognition of their immense potential, and their limitations. Harnessing their synergy will lead to breakthroughs and spectacular progress. This chapter is a modest attempt to lay some ground for this symbiosis.

\subsection{Literature review}
Economics and finance are two of the fields of choice to which the early contributors to the theory of Mean Field Games (MFGs) and Mean Field Control (MFC) paid special attention. To wit, more than $30$ pages of the introductory chapter \emph{"Learning by Examples"} of \cite{MR3752669} were devoted to applications to these fields. Additionally, a more recent review of the literature which appeared on this very subject since 2018 can be found in \cite{carmona-AMS-tmp}.
Here we shall mention some of the most striking applications, and emphasize those which triggered progress in numerical methods and especially applications of Machine Learning (ML) techniques.
\vskip 2pt
Historically, macro-economic models have been cast as general equilibrium problems and solved as such. However, many of these models, see for example \cite{krusell-smith,aiyagari,bewley,huggett}, carry all the elements of MFGs, and have since been revisited in light of our understanding of MFGs. See for example \cite{MR3679338,nuno2018social} and \cite{sargent2020stochastic}, and \cite{MR3268061,achdou2017income} for a numerical point of view.
\vskip 2pt
But as emphasized in \cite{carmona-AMS-tmp}, the crusade of Brunnermeier and Sannikov arguing for the merging of macro-economics and finance models through the common use of continuous time helped the convergence of economists, financial engineers and applied mathematicians toward the use of a common language and a common set of models, MFG models playing a crucial role in this evolution.
\vskip 2pt
The analysis of systemic risk inherent in large banking networks was a natural ground for mean field models of interactions. The early work \cite{MR3032938} and the more recent model \cite{MR3910001} lead to challenging Partial Differential Equation (PDE) problems, while \cite{MR3325083} offers a simple model which can be solved explicitly, both in its finitely many player version and its infinite player form, and for which the master equation can be derived and solved explicilty. Subsequent and more realistic models involving delays like \cite{MR3865333} or interactions through default times like~\cite{ElieIchibaLauriere-2020-LargeBankingMFG} are unfortunately more difficult to solve, even numerically.
\vskip 2pt
Some of the applications of mean field games lie at the intersection of macro-economics and  financial engineering. As an example, we mention the MFG models for Bertrand and Cournot equilibria of oil production introduced in \cite{Gueant-2013-phdthesis,MR2762362}, 
and revisited later on in~\cite{MR3359708,MR3683681} with an interesting emphasis on exhaustability, and a more mathematical PDE analysis of the model in \cite{MR3755719}.
\vskip 2pt
Like the macro-economic general equilibrium models mentioned earlier, models of bank run such as \cite{rochet2004coordination} are \emph{screaming} for a MFG reformulation and this was first done in \cite{MR3780736} and \cite{MR3679343}, and later on in \cite{MR3906158} by analytic methods, introducing a new class of MFGs of timing. 
However, like in macro-economics realistic models require the introduction of a common noise, making the theoretical solution much more involved, and numerical implementations even more difficult.
\vskip 2pt
The high frequency markets offer, without the shadow of a doubt, one of the most computer intensive financial engineering application one can think of. While most existing papers on the subject revolve around the properties of price impact, see e.g. \cite{MR3325272,MR3805247,MR3500455}, modeling the interaction between a large number of market makers and trading programs is certainly a exciting challenge that the methodology and the numerical tools developed for the analysis of MFGs should make it possible to meet.

	\vskip 2pt
	The introduction of the MFG paradigm opened the door to the search for solutions of large population equilibrium models which could not be imagined to be solvable before. Still, the actual solutions of practical applications had to depend on the development of efficient numerical algorithms implementing the MFG principles. This was done early in the development of the theory. See e.g.~\cite{MR2679575,achdoulauriere-2020-mfg-numerical,AMSnotesLauriere} and the references therein. More recently, the release in the public domain of powerful software packages such as {\tt TensorFlow} has made it possible to test at a very low cost the possible impact of machine learning (ML) tools in the solution of challenging problems for MFGs and MFC, whether these problems were formulated in the probabilistic approach~\cite{fouque2019deep,carmona2019convergence2,germain2019numerical} or the analytical approach~\cite{al2018solving,carmona2021convergence1,ruthotto2020machine,cao2020connecting,lin2020apac,AMSnotesLauriere}. These new methods combine neural network approximations and stochastic optimization techniques to solve McKean-Vlasov control problems, mean field forward-backward stochastic differential equations (FBSDE) or mean field PDE systems. 
	
	\vskip 2pt
	While the present chapter concentrates on ML applications of MFGs and MFC in finance, the reader should not be surprised if they recognize a strong commonality of ideas and threads with the subchapter~\cite{angiuli2021reinforcement} dealing with reinforcement learning for MFGs with a special focus on a two-timescale procedure, and the subchapter~\cite{germain2021neural} offering a review of neural-network-based algorithms for stochastic control and PDE applications in finance.

\vskip 2pt
The rest of this subchapter is organized as follows. In the rest of this section, we define the MFG and MFC problems. In Section~\ref{sec:directMethod}, we present a direct ML method for MKV control. We then turn our attention to method related to the optimality conditions of MFGs and MFC. A neural-network based shooting method for generic MKV FBSDE systems is discussed in Section~\ref{sec:MFdeepBSDE}. A deep learning method for mean-field PDE systems is presented in Section~\ref{sec:MFDGM}. We conclude in Section~\ref{sec:ccl}.

\subsection{Definition of the problems}

The parameters of our models are a time horizon $T>0$, integers $d$ and $k$ for the dimensions of the state space 
$\dom \subseteq \RR^d$ and the action space $\RR^k$. Typically, the state space $\dom$ will be the whole space $\RR^d$. We shall use the notation $\domT = [0,T] \times \dom$, and $\langle \cdot, \cdot \rangle$ for the inner product of two vectors of compatible dimensions.
We denote by $\cP_2(\dom)$ the space of probability measures on $\dom$ which integrate the square of the norm when $\dom$ is unbounded.

Now let $f: \dom \times \cP_2(\dom) \times \RR^\ctrldim \to \RR, (x,m,\ctrl) \mapsto f(x,m,\ctrl)$ and $g: \dom \times \cP_2(\dom) \to \RR, (x,m) \mapsto g(x,m)$ be functions giving respectively the instantaneous running cost and the terminal cost, let $b:  \dom \times \cP_2(\dom) \times \RR^\ctrldim \to \RR^d, (x,m,\ctrl) \mapsto b(x,m,\ctrl)$ be its drift function, and let $\sigma>0$ be the volatility of the state's evolution (for simplicity we focus on the case of a non-degenerate diffusion although some of the methods presented below can also be applied when the diffusion is degenerate).  These functions could be allowed to also depend on time at the expense of heavier notation. Here, $x,m$ and $\ctrl$ play respectively the role of the state of the agent, the mean-field term (i.e. the population's distribution), and the control used by the agent. In general, the mean-field term is a probability measure. However, in some cases, we will assume that this probability measure has a density which is in $L^2(\dom)$.
\vskip 2pt

\begin{definition}[MFG equilibrium] 
\label{def:MFGeq}
When considering the \index{mean field game}mean field game problem for a given initial distribution $m_0 \in \cP_2(\dom)$, we call a Nash equilibrium a flow of probability measures $\hat{m}=(\hat m(t,\cdot))_{0\le t\le T}$ in $\cP_2(\dom)$ and a feedback control $\hat{\ctrl}: \domT\to \RR^\ctrldim$ satisfying the following two conditions:
\begin{enumerate}
	\item $\hat{\ctrl}$ minimizes $J^{MFG}_{\hat{m}}$ where, for $m=(m(t,\cdot))_{0\le t\le T}$,
\begin{align}
\label{subchapRCML-num-eq:def-J-MFG}
	J^{MFG}_{m}: \ctrl \mapsto  \EE \left[\int_0^T f(X_t^{m, \ctrl}, m(t, \cdot), \ctrl(t,X_t^{m, \ctrl}) ) dt + g(X_T^{m, \ctrl}, m(T,\cdot)) \right]
\end{align}
under the constraint that the  process $X^{m, \ctrl} = (X_t^{m, \ctrl})_{t \ge 0}$ solves the stochastic differential equation (SDE)
\begin{equation}
\label{subchapRCML-num-eq:dyn-X-general-MFG}
	d X_t^{m, \ctrl} = b(X_t^{m, \ctrl}, m(t,\cdot), \ctrl(t, X_t^{m, \ctrl})) dt + \sigma d W_t, \qquad t \ge 0,
\end{equation}
where $W$ is a standard $d$-dimensional Brownian motion, and $X_0^{m, \ctrl}$ has distribution $m_0$;
	\item For all $t \in [0,T]$, $\hat{m}(t,\cdot)$ is the probability distribution of $X_t^{\hat{m}, \hat{\ctrl}}$.
\end{enumerate}
\end{definition}
In the definition~\eqref{subchapRCML-num-eq:def-J-MFG} of the cost function, the subscript $m$ is used to emphasize the dependence on the mean-field flow, which is fixed when an infinitesimal agent performs their optimization. The second condition ensures that if all the players use the control $\hat{\ctrl}$ identified in the first bullet point, the law of their individual states is indeed $\hat{m}$.

\vskip 2pt
Using the same drift, running and terminal cost functions and volatility, we can also consider the corresponding McKean-Vlasov (MKV for short) control, or mean field control (MFC for short) problem. This optimization problem corresponds to a social optimum and is phrased as an optimal control problem. It can be interpreted as a situation in which all the agents cooperate to minimize the average cost. 

\begin{definition}[MFC optimum]
A feedback control $\ctrl^*: \domT\to \RR^\ctrldim$ is an optimal control for the MKV control (or MFC)  problem for a given initial distribution $m_0 \in \cP_2(\dom)$ if it minimizes
\begin{align}
\label{subchapRCML-num-eq:def-J-MFC}
	J^{MFC}: \ctrl \mapsto \EE \left[\int_0^T f(X_t^{\ctrl}, m^{\ctrl}(t,\cdot), \ctrl(t,X_t^{\ctrl}) ) dt + g(X_T^{\ctrl}, m^{\ctrl}(T,\cdot)) \right]
\end{align}
where $m^{\ctrl}(t,\cdot)$ is the probability distribution of the law of $X_t^{\ctrl}$, under the constraint that the process $X^{\ctrl} = (X_t^{\ctrl})_{t \ge 0}$ solves the stochastic differential equation of McKean-VLasov type:
\begin{equation}
\label{subchapRCML-num-eq:dyn-X-general-MFC}
	d X_t^{\ctrl} = b(X_t^{\ctrl}, m^{\ctrl}(t,\cdot), \ctrl(t, X_t^{\ctrl})) dt + \sigma d W_t, \qquad t \ge 0,
\end{equation}
$X_0^{\ctrl}$ having distribution $m_0$. 
\end{definition}
The application of MFC are not limited to social optima in very large games. These problems also arise for example in risk management~\cite{MR2784835} or in optimal control with a cost involving a conditional expectation~\cite{achdoulaurierelions2020optimal,MR4133380}.
\vskip 0pt
If $m^* = m^{\ctrl^*}$ is the flow of state distribution for an optimal control $\ctrl^*$, then:
$$
	J^{MFC}(\ctrl^*) = J^{MFG}_{m^*}(\ctrl^*) \leq J^{MFG}_{\hat{m}}(\hat{\ctrl}).
$$
In general the inequality is strict, which leads to the notion of price of anarchy~\cite{MR3968548}. 

To simplify the presentation, we have introduced MFG and MFC in a basic formulation where the interactions occur through the distribution of states. However, in many applications, the interactions occur through the distribution of controls or through the joint distribution of states and controls. This aspect will be illustrated in some of the examples discussed below.

\section{Direct method for MKV Control}
\label{sec:directMethod}

In this section, we present a direct approach to the numerical solution of McKean-Vlasov control problems. It hinges on an approach developed for standard control problems, in which the control feedback function is 
restricted to a parametric family of function, especially a class of neural networks whose parameters are learned by stochastic optimization~ \cite{MR2137498,han2016deep-googlecitations}.  This method was extended to the mean-field setting in~ \cite{fouque2019deep,carmona2019convergence2}. We illustrate this approach with the solution of a price impact model.

\subsection{Description of the method}

Since MFC is an optimization problem, it is natural to leverage stochastic optimization tools from machine learning directly applied to the definition~\eqref{subchapRCML-num-eq:def-J-MFC}--\eqref{subchapRCML-num-eq:dyn-X-general-MFC}. We introduce three approximations leading to a formulation more amenable to numerical treatment. 

\vskip 10pt
\textbf{Approximation steps.} First, we restrict the set of (feedback) controls to be the set of neural networks with a given architecture. We introduce new notation to define this class of controls. 
We denote by:
\begin{align*}
	\mathbf{L}^\psi_{d_1, d_2} = 
	&\Big\{ \phi: \RR^{d_1} \to \RR^{d_2} \,\Big|\,  \exists (\beta, w) \in  \RR^{d_2} \times \RR^{d_2 \times d_1}, \forall  i \in \{1,\dots,d_2\}, 
	\\
	&\qquad \;\phi(x)_i = \psi\Big(\beta_i + \sum_{j=1}^{d_1} w_{i,j} x_j\Big) \Big\} 
\end{align*}
the set of layer functions with input dimension $d_1$, output dimension $d_2$, and activation function $\psi: \RR \to \RR$. Typical choices for $\psi$ are the ReLU function $\psi(x) = x^+$ or the sigmoid function $\psi(x) = 1/(1+e^{-x})$. 
Building on this notation and denoting by $\circ$ the composition of functions, we define:
\begin{align*}
	\bN^\psi_{d_0, \dots, d_{\ell+1}} 
	= 
	&\Big\{ \varphi: \RR^{d_0} \to \RR^{d_{\ell+1}} \,\Big|\, \exists (\phi_i)_{i=0, \dots, \ell-1} \in \bigtimes_{i=0}^{i=\ell-1} \mathbf{L}^\psi_{d_i, d_{i+1}}, 
	\\
	&\qquad \exists \phi_\ell \in \mathbf{L}^{\mathrm{id}}_{d_{\ell}, d_{\ell+1}}, \varphi = \phi_\ell \circ \phi_{\ell-1} \circ \dots \circ \phi_0  \Big\} \, 
	\notag
\end{align*}
 the set of regression neural networks with $\ell$ hidden layers and one output layer, the activation function of the output layer being the identity $\psi=\text{id}$. The number $\ell$ of hidden layers, the numbers $d_0$, $d_1$, $\cdots$ , $d_{\ell+1}$ of units per layer, and the activation functions, is what is called the architecture of the network. Once it is fixed, the actual network function $\varphi\in \bN^\psi_{d_0, \dots, d_{\ell+1}} $ is determined by the remaining real-valued parameters:
 $$
 \theta=(\beta^{(0)}, w^{(0)},\beta^{(1)}, w^{(1)},\cdots\cdots,\beta^{(\ell-1)}, w^{(\ell-1)},\beta^{(\ell)}, w^{(\ell)})
 $$
defining the functions $\phi_0$, $\phi_1$, $\cdots$ , $\phi_{\ell-1}$ and $\phi_\ell$ respectively. The set of such parameters is denoted by $\Theta$. For each $\theta\in\Theta$, the function $\varphi$ computed by the network will be denoted by $\ctrl_\theta \in \bN^\psi_{d_0, \dots, d_{\ell+1}}$.
As it should be clear from the discussion of the previous section, here, we are interested in the  case where $d_0 = d+1$ (since the inputs are time and state) and $d_{\ell+1} = \ctrldim$ (i.e., the control dimension). 

\vskip 2pt
Our first approximation is to minimize $J^{MFC}$ defined by~\eqref{subchapRCML-num-eq:def-J-MFC}--\eqref{subchapRCML-num-eq:dyn-X-general-MFC} over $\ctrl \in \bN^\psi_{d+1, d_1, \dots, d_{\ell}, \ctrldim} $, or equivalently, to minimize over $\theta \in \Theta$ the function:
\begin{align*}
	\bJ: \theta \mapsto \EE \left[\int_0^T f(X_t^{\ctrl_\theta}, m^{\ctrl_\theta}(t,\cdot), \ctrl_\theta(t,X_t^{\ctrl_\theta}) ) dt + g(X_T^{\ctrl_\theta}, m^{\ctrl_\theta}(T,\cdot)) \right]
\end{align*}
where $m^{\ctrl_\theta}(\cdot, t)$ is the  law of $X_t^{\ctrl_\theta}$, under the constraint that the process $X^{\ctrl_\theta}$ solves the SDE~\eqref{subchapRCML-num-eq:dyn-X-general-MFC} with feedback control $\ctrl_\theta$.

\vskip 2pt
Next, we approximate the probability distribution of the state. A (computationally) simple option is to replace it by the empirical distribution of a system of $N$ interacting particles. Given a feedback control $\ctrl$, we denote by $(\underline X_t^{\ctrl})_{t} = (X_t^{1, \ctrl}, \dots, X_t^{N, \ctrl})_{t}$ the solution of the system:
\begin{equation}
\label{subchapRCML-num-eq:NN-MFC-Nparticles}
	d X_t^{i, \ctrl} = b(X_t^{i, \ctrl}, m^{N, \ctrl}_t, \ctrl(t, X_t^{i, \ctrl})) dt + \sigma d W^i_t, \qquad t \ge 0,\;\;i=1,\ldots,N
\end{equation}
where 
$$
	m^{N, \ctrl}_t = \frac{1}{N} \sum_{j=1}^N \delta_{X_t^{j, \ctrl}}, 
$$ 
is the empirical measure of the $N$ particles, $(W^i)_{i=1,\dots,N}$ is a family of $N$ independent $d$-dimensional Brownian motions, and the initial positions $(X_0^{i, \ctrl})_{i=1,\dots,N}$ are i.i.d. with distributions $m_0$. The $N$ stochastic differential equations in~\eqref{subchapRCML-num-eq:NN-MFC-Nparticles} are coupled via their drifts through the empirical measure $m^{N, \ctrl}_t$. The controls are distributed in the sense that the control used in the equation for $X^{i, \ctrl}$ is a function of $t$ and $X_t^{i, \ctrl}$ itself, and not of the states of the other particles. Despite their dependence due to the coupling, it is expected that the empirical measures converge when $N\to\infty$ to the solution of the SDE~\eqref{subchapRCML-num-eq:dyn-X-general-MFC}. Not only does this convergence holds, but in this limit, the individual particle processes $( X_t^{i, \ctrl})_{0\le t\le T}$ become independent in this limit. This fundamental result is known under the name of \emph{propagation of chaos}. See~\cite{MR3753660} for details and the role this result plays in the theory of MFGs and MFC. As per this second approximation, the new problem is to minimize over $\theta \in \Theta$ the function
\begin{align*}
	\bJ^{N}: \theta \mapsto \frac{1}{N} \sum_{i=1}^N \EE \left[\int_0^T f(X_t^{i, \ctrl_\theta}, m^{N, \ctrl_\theta}_t, \ctrl_\theta(t,X_t^{i, \ctrl_\theta}) ) dt + g(X_T^{i, \ctrl_\theta}, m^{N,\ctrl_\theta}_T) \right],
\end{align*}
under the dynamics~\eqref{subchapRCML-num-eq:NN-MFC-Nparticles} with control $\ctrl_\theta$.
\vskip 2pt
Our third approximation is to discretize time. Let $N_T$ be a positive integer, let $\Delta t = T/N_T$ and $t_n = n \Delta t$, $n =0,\dots, N_T$. We now minimize over $\theta \in \Theta$ the function
\begin{equation}
\label{subchapRCML-num-eq:NN-MFC-cost-totalapprox}
	\bJ^{N, \Delta t}: \theta \mapsto \frac{1}{N} \sum_{i=1}^N \EE \left[\sum_{n=0}^{N_T-1} f(\check X_{t_n}^{i, \ctrl_\theta}, \check m^{N, \ctrl_\theta}_{t_n}, \ctrl_\theta(t_n,\check X_{t_n}^{i, \ctrl_\theta}) ) \Delta t + g(\check X_T^{i, \ctrl_\theta}, \check m^{N,\ctrl_\theta}_T) \right],
\end{equation}
under the dynamic constraint:
\begin{equation}
\label{subchapRCML-num-eq:NN-MFC-Nparticles-Deltat}
	\check X_{t_{n+1}}^{i, \ctrl_\theta} = \check X_{t_{n}}^{i, \ctrl_\theta} + b( \check X_{t_{n}}^{i, \ctrl_\theta}, \check m^{N, \ctrl_\theta}_{t_{n}}, \ctrl_\theta(t, \check X_{t_{n}}^{i, \ctrl_\theta})) \Delta t + \sigma \Delta \check W^i_n, \qquad n = 0, \dots, N_T-1,
\end{equation}
and the initial positions $(\check X_0^{i, \ctrl_\theta})_{i=1,\dots,N}$ are i.i.d. with distribution given by the density $m_0$, where 
$$
	\check m^{N, \ctrl_\theta}_{t_{n}} = \frac{1}{N} \sum_{j=1}^N \delta_{\check X_{t_{n}}^{j, \ctrl_\theta}}, 
$$ 
and the $(\Delta \check W_n^i)_{i=1,\dots,N,n=0,\dots,N_T-1}$ are i.i.d. random variables with Gaussian distribution $\mathcal N(0, \Delta t)$.

Under suitable assumptions on the model and the neural network architecture, the difference between $\inf_\theta\bJ^{N, \Delta t}(\theta)$ and $\inf_\ctrl J^{MFC(\ctrl)}$ goes to $0$ as $N_t, N$ and the number of parameters in the neural network go to infinity. See~\cite{carmona2019convergence2} for details.

\vskip 10pt
\textbf{Optimization procedure. } Two obvious difficulties have to be overcome. First the fact that the cost function~\eqref{subchapRCML-num-eq:NN-MFC-cost-totalapprox} is in general non-convex. Second, the parameter $\theta$ is typically high dimensional. But the cost~\eqref{subchapRCML-num-eq:NN-MFC-cost-totalapprox} being written as an expectation, it is reasonable to rely on a form of stochastic gradient descent (SGD) algorithm. The randomness in our problem comes from the initial positions $\underline {\check X}_0 = (\check X_0^{i, \ctrl_\theta})_{i}$ and the random shock innovations $(\Delta \underline {\check W}_n)_{n=0,\dots,N_T} = (\Delta \check W_n^i)_{i,n}$. Hence $S = (\underline {\check X}_0, (\Delta \underline {\check W}_n)_{n})$ is going to play the role of a random sample in SGD. Given  a realization of $S$ and a choice of parameter $\theta$, we can construct the trajectory $(\check X_{t_n}^{i, \ctrl_\theta,S})_{i=1,\dots,N, n=0,\dots,N_T}$ by following~\eqref{subchapRCML-num-eq:NN-MFC-Nparticles-Deltat}, and compute the induced cost:
\begin{equation}
\label{subchapRCML-num-eq:NN-MFC-cost-totalapprox-oneS}
	\bJ^{N, \Delta t}_S(\theta) = \frac{1}{N} \sum_{i=1}^N \left[\sum_{n=0}^{N_T-1} f(\check X_{t_n}^{i, \ctrl_\theta, S}, \check m^{N, \ctrl_\theta, S}_{t_n}, \ctrl_\theta(t_n,\check X_{t_n}^{i, \ctrl_\theta, S}) ) \Delta t + g(\check X_T^{i, \ctrl_\theta, S}, \check m^{N,\ctrl_\theta, S}_T) \right].
\end{equation}
We give our SGD procedure in Algorithm~\ref{subchapRCML-num-algo:SGD-MFC}. The most involved step is the computation of the gradient $\nabla_\theta \bJ^{N, \Delta t}_S(\theta^{(\mathtt{k})})$  with respect to $\theta$. However, modern programming libraries (such as {\tt TensorFlow} or {\tt PyTorch}) perform this computation automatically using backpropagation, simplifying dramatically the code. The present method is thus extremely straightforward to implement: contrary to the methods based on optimality conditions there is no need to derive by hand any PDE, any FBSDE, or compute gradients. We work directly with the definition of the MFC.

Besides this aspect, the main reasons behind the success of this method are the expressive power of neural networks and the fact that there is a priori no limitation on the number $\mathtt{K}$ of iterations because the samples $S$ come from Monte-Carlo simulations and not from a training set of data. In the implementation of this method, using mini-batches and ADAM~\cite{kingma2014adam} can help improving convergence.

\vskip 6pt
\begin{algorithm}[ht]
\DontPrintSemicolon
\KwData{An initial parameter $\theta^{(0)}\in\Theta$; a number of steps $\mathtt{K}$;  a sequence $(\beta^{(\mathtt{k})})_{\mathtt{k}=0,\dots,\mathtt{K}-1}$ of learning rates.}
\KwResult{A parameter $\theta$ such that $\ctrl_{\theta}$ approximately minimizes $J^{MFC}$}
\Begin{
  \For{$\mathtt{k} = 0, 1, 2, \dots, \mathtt{K}-1 $}{
    Pick $S  = (\underline {\check X}_0, (\Delta \underline {\check W}_n)_{n})$\;
    Compute the gradient $\nabla_\theta \bJ^{N, \Delta t}_S(\theta^{(\mathtt{k})})$, see~\eqref{subchapRCML-num-eq:NN-MFC-cost-totalapprox-oneS}\;
    Set $\theta^{(\mathtt{k}+1)} = \theta^{(\mathtt{k})} - \beta^{(\mathtt{k})} \nabla_\theta \bJ^{N, \Delta t}_S(\theta^{(\mathtt{k})})$ \;
    }
  \KwRet{ $\theta^{(\mathtt{K})}$}
  }
\caption{SGD for MFC\label{subchapRCML-num-algo:SGD-MFC}}
\end{algorithm}

\subsection{Numerical illustration: a price impact model}
\label{sub:price_impact}
We consider a financial application originally solved as a mean field game by Carmona and Lacker in the weak formulation in \cite{MR3325272}, and revisited in the book of Carmona and Delarue  \cite[Sections 1.3.2 and 4.7.1]{MR3752669} in the strong formulation. This is a problem of optimal execution in the presence of price impact resulting from a large group of traders affecting the price of a single asset through their aggregate trades: a large number of buy orders will push the price up while a large number of sell orders will deflate the price.
 This aggregate impact on the price is the source of the \emph{mean field} in the model. As a consequence, this model is an instance of mean field problems with interactions through the distribution of controls, introduced by Gomes et al. in~\cite{MR3160525,MR3489048} who coined the term \emph{extended MFG}. 
 
 Here we shall consider the MFC counterpart. The method described above can readily be adapted to solve MFC with interactions through the control's distribution by computing the empirical distribution of controls for an interacting system of $N$ particles. For the sake of completeness, we recall the model and the derivation of a succinct MFC formulation. See the aforementioned references for the $N$-agent problem and more details (in the MFG setting). A typical and infinitesimal trader's inventory at time $t$ is denoted by $X_t$.  We assume that it evolves according to the SDE:
\begin{equation*}
    dX_t = \ctrl_t dt +\sigma dW_t,
\end{equation*}
where $\ctrl_t$ represents the rate of trading, and $W$ is a standard Brownian motion.  

Denoting by $\nu^\ctrl_t=\mathcal{L}(\ctrl_t)$ the law of the control at time $t \in [0,T]$, we consider that the evolution of the price is given by:
\begin{equation*}
    dS_t =\gamma \biggl( \int_{\mathbb{R}}a d\nu^\ctrl_t(a) \biggr) dt + \sigma_0 dW_t^0, 
\end{equation*} 
where $\gamma$ and $\sigma_0$ are positive constants, and the Brownian motion $W^0$ is independent from $W$. The \emph{price impact} effect is taken into account by the fact that the price $(S_t)_{0 \leq t \leq T}$ of the asset is influenced by the average trading rate of all the traders.

If we denote by $K_t$ the amount of cash held by the trader at time $t$, the dynamics of the process $(K_t)_{0 \le t \le T}$ are given by:
\begin{equation*}
    dK_t=-\big(\ctrl_t S_t +c_{\ctrl}(\ctrl_t)\big)dt,
\end{equation*}
where the function $a \mapsto c_{\ctrl}(a)$ is a non-negative convex function satisfying $c_{\ctrl}(0)=0$. It corresponds to the cost for trading at rate $a$. At time $t$, the trader's total wealth, denoted by $V_t$, is the sum of the cash holding and the value of the inventory marked at the price $S_t$, i.e., $V_t=K_t+X_t S_t$. 
Using the self-financing condition of Black-Scholes' theory, the evolution of the trader's wealth is:
\begin{equation}
\begin{split}
    dV_t&=dK_t+X_tdS_t+S_t dX_t
    \\
   & =\Big[ -c_{\ctrl}(\ctrl_t)+\gamma X_t \int_{\mathbb{R}} a d\nu^\ctrl_t(a) \Big]dt + \sigma S_t dW_t + \sigma_0 X_t dW_t^0.
\end{split}
\label{eq:wealth}
\end{equation}
We assume that the trader is subject to a running cost for holding an inventory, modeled by a function $c_X$ of their inventory, and to a terminal liquidation constraint at maturity $T$ represented by a scalar function $g$. Thus, the trader's cost function, to be minimized, is defined by:
\begin{equation*}
J(\ctrl)=\mathbb{E}\Big[\int_0^T c_X(X_t) dt +g(X_T) - V_T\Big].
\end{equation*}
Taking into account~\eqref{eq:wealth}, this cost can be rewritten in terms of $X$ only as:
\begin{equation*}
J(\ctrl)=\mathbb{E}\Big[ \int_0^T \left(c_{\ctrl}(\ctrl_t)+c_X(X_t)-\gamma X_t \int_{\mathbb{R}} a  d\nu^\ctrl_t(a)\right)dt +g(X_T)\Big].  
\end{equation*}
Following the Almgren-Chriss linear price impact model, we assume that the functions $c_X$, $c_{\ctrl}$ and $g$ are quadratic. Thus, the cost is of the form:
\begin{equation*}
J(\alpha)=\mathbb{E}\left[ \int_0^T \left( \frac{c_{\alpha}}{2}{\alpha_t}^2+\frac{c_X}{2}X_t^2-\gamma X_t\int_{\mathbb{R}} a d\nu^\ctrl_t(a) \right)dt + \frac{c_g}{2}X_T^2\right].
\end{equation*}
Let us stress that this problem is an extended MFC: the population distribution is not frozen during the optimization over $\ctrl$, and the interactions occur through the distribution of controls $\nu^\ctrl$.

\begin{figure}[ht]
\centering
\begin{subfigure}{.45\textwidth}
  \centering
  \includegraphics[width=\linewidth]{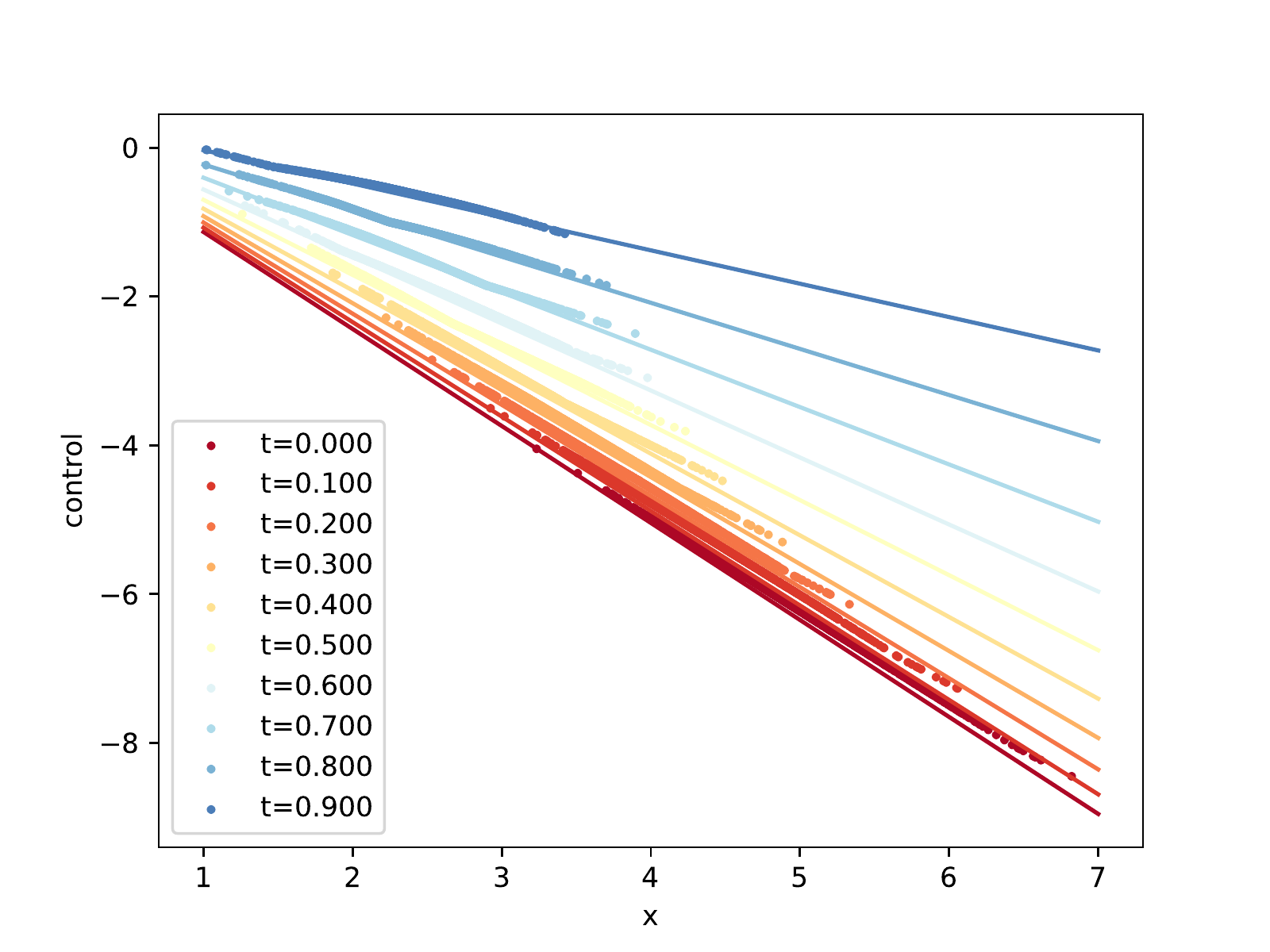}
\end{subfigure}
\begin{subfigure}{.45\textwidth}
  \centering
  \includegraphics[width=\linewidth]{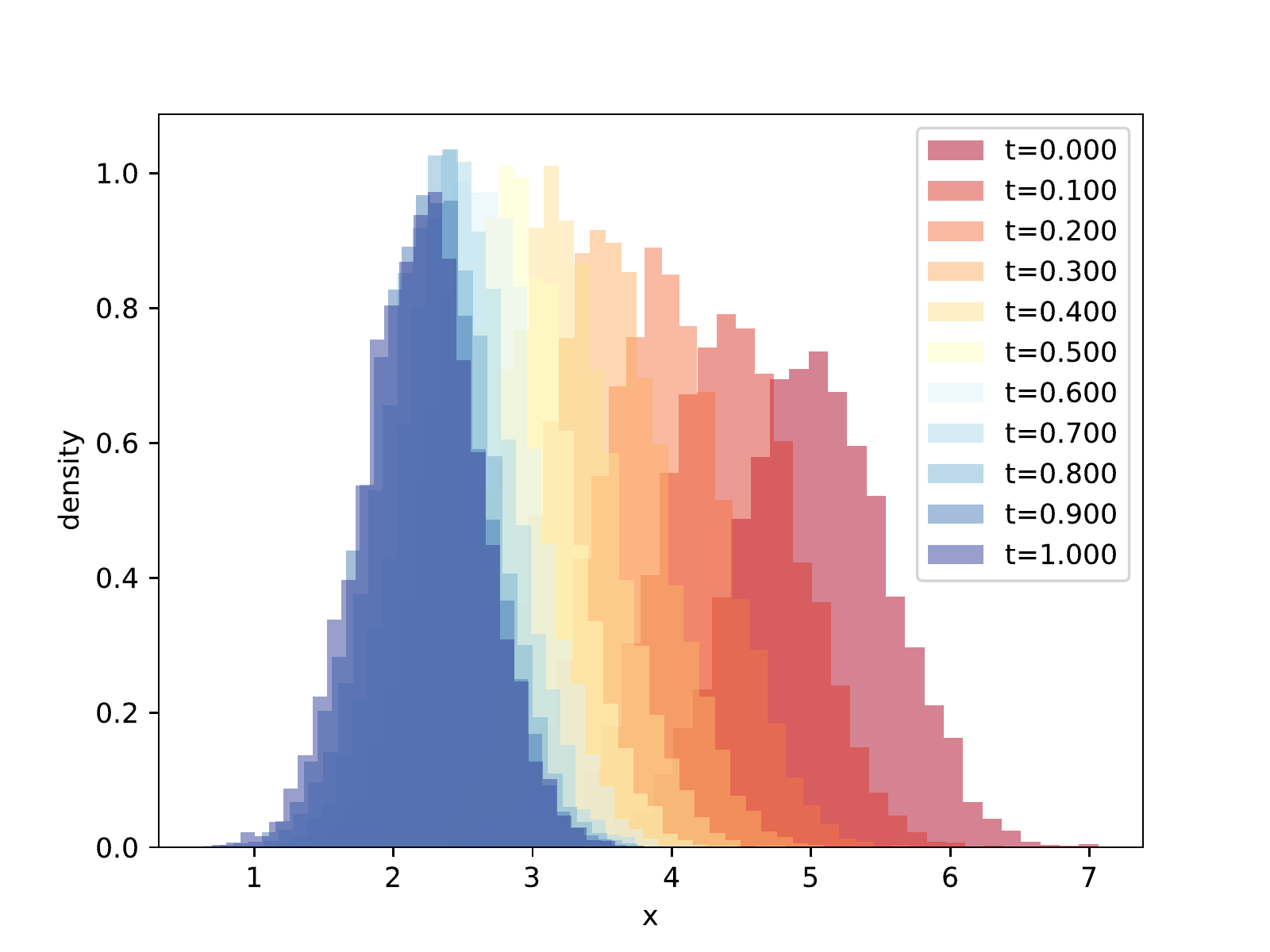}
\end{subfigure}
\caption{Price impact MFC example solved by Algorithm~\ref{subchapRCML-num-algo:SGD-MFC}. Left: Control learnt (dots) and exact solution (lines). Right: associated empirical state distribution. Here, $\gamma = 0.2$.}
\label{fig:ex-price-impact-1}
\end{figure}

\begin{figure}[ht]
\centering
\begin{subfigure}{.45\textwidth}
  \centering
  \includegraphics[width=\linewidth]{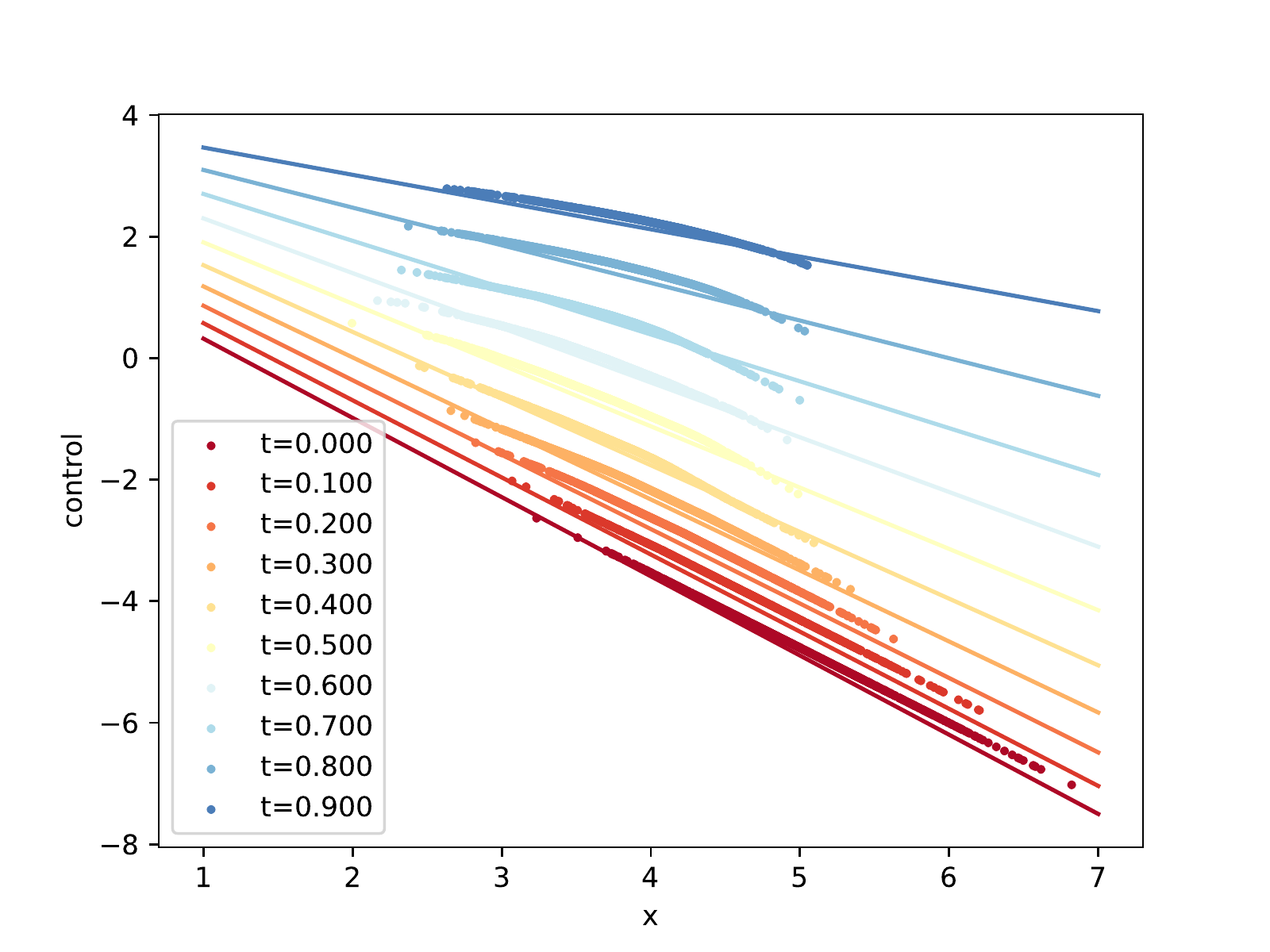}
\end{subfigure}
\begin{subfigure}{.45\textwidth}
  \centering
  \includegraphics[width=\linewidth]{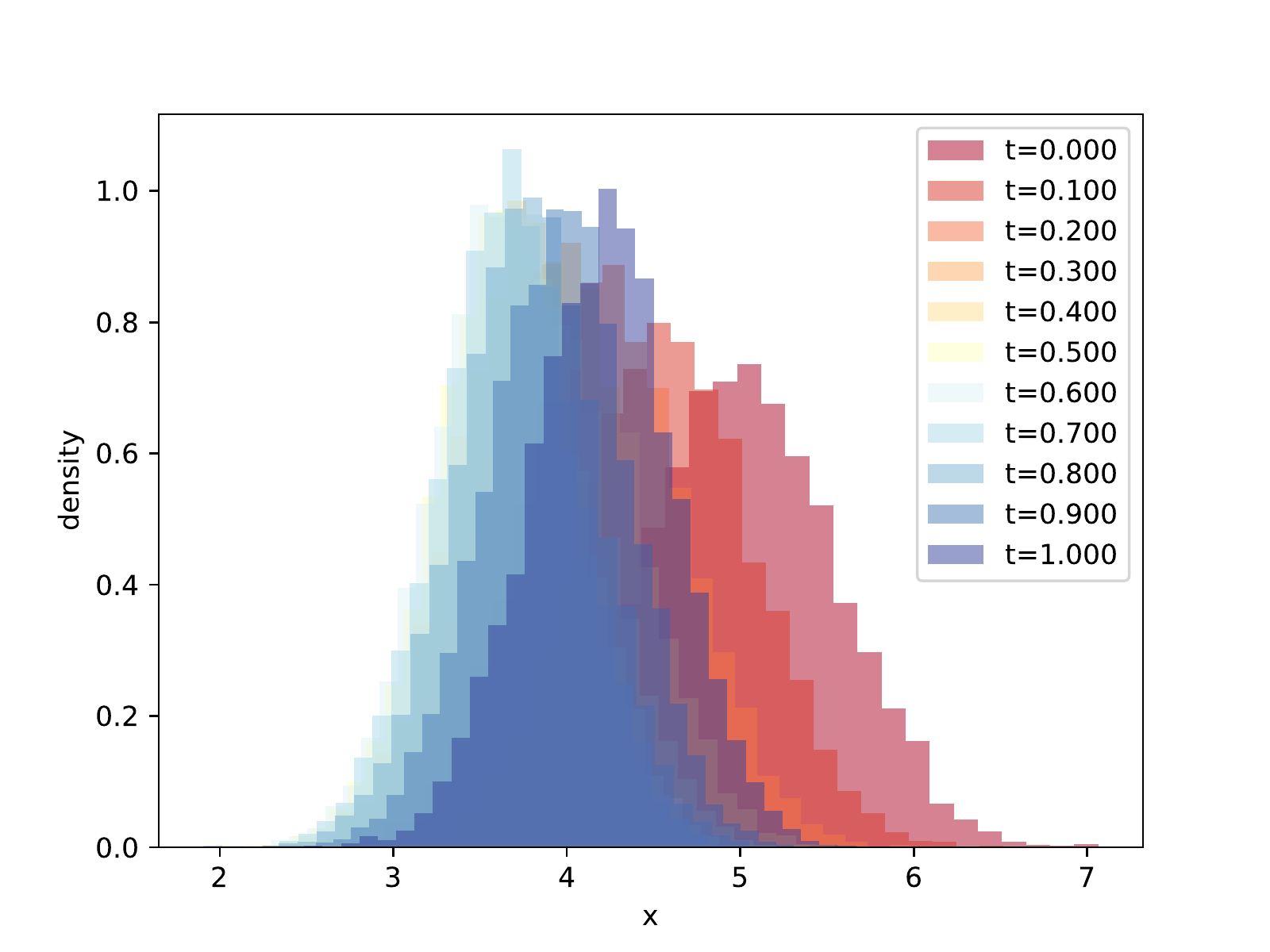}
\end{subfigure}
\caption{Price impact MFC example solved by Algorithm~\ref{subchapRCML-num-algo:SGD-MFC}. Left: Control learnt (dots) and exact solution (lines). Right: associated empirical state distribution. Here,  $\gamma = 1$.}
\label{fig:ex-price-impact-2}
\end{figure}

This model is solved by reinforcement learning techniques (for both MFC and the corresponding MFG) in~\cite{angiuli2021reinforcement}. Here, we present results obtained with the deep learning method described above, see Algorithm~\ref{subchapRCML-num-algo:SGD-MFC}. The results shown in Figure~\ref{fig:ex-mfg-sysrisk-traj} show that the control is linear, as expected from the theory, and the distribution moves towards $0$ while becoming more concentrated. In other words, at the beginning the traders have a relatively large inventory with a large variance across the population, and they liquidate to end up with smaller inventories and less variance. One can see that towards the end of the time interval, the control learnt is not exactly linear, probably because a regions has been less explored than the rest leading to a less accurate training. For these results, we used the parameters: $T=1$, $c_X = 2$, $c_{\ctrl} = 1$, $c_g = 0.3$, $\sigma = 0.5$ and the value of $\gamma$ indicated in the captions. 
Moreover, in the algorithm we took $N=2000$ particles and $N_T = 50$ time steps. We see in Figure~\ref{fig:ex-price-impact-1} that when $\gamma=0.2$, the optimal control is to constantly liquidate. However, as shown in Figure~\ref{fig:ex-price-impact-2}, when $\gamma=1$, the traders start by liquidating by towards the end of the time interval, they buy. This can be explained by the fact that with a higher $\gamma$, the price impact effect is stronger and the traders can use phenomenon to increase their wealth by collectively buying and hence increasing the price. In each case, the neural network manages to learn a control which approximately matches the semi-explicit one obtained by reducing the problem to an system of ordinary differential equations (ODE) as explained in~\cite[Sections 1.3.2 and 4.7.1]{MR3752669}.

\section{Deep BSDE method for MKV FBSDEs}
\label{sec:MFdeepBSDE}

In this section, we present an extension to the mean-field regime of the DeepBSDE method introduced in~\cite{MR3736669} and analyzed in~\cite{han2020convergence}. The latter uses neural networks to learn the solution of BSDEs. It relies on a shooting method, where one tries to find a suitable starting point in order to match a given terminal condition. These ideas have been extended to the mean-field setting to solve forward-backward systems of McKean-Vlasov SDEs in~\cite{fouque2019deep,carmona2019convergence2,germain2019numerical}. After presenting the main ideas, we illustrate the performance of the method on a systemic risk MFG model introduced in~\cite{MR3325083} for which explicit solutions exist.

\subsection{Description of the method}

As explained in \cite{MR3752669}, finding an MFG equilibrium can be reduced to the solution of a forward-backward system of SDEs (FBSDE for short) which reads:
\begin{equation}
\label{eq:MKV-FBSDE}
\left\{
\begin{aligned}
	d X_t
	=
	& b\bigl(t,X_t,\cL(X_t),\hat{\alpha}(t, X_t, \cL(X_t), Y_t)\bigr) dt \,+\sigma dW_t
	\\
	d Y_t
	=
	& - \partial_x H(t, X_t, \cL(X_t), Y_t,Z_t,\hat{\alpha}(t, X_t, \cL(X_t), Y_t))  dt
+ Z_t d W_t,
\end{aligned}
\right.
\end{equation}
with initial condition $X_0$ having distribution $m_0$ 
 and terminal condition $Y_T = \partial_x g(X_T, \cL(X_T))$. 
 $H$ is the Hamiltonian:
$$
H(t,x,\mu,y,z,\alpha)=b(t,x,\mu,\alpha)\cdot y + \sigma\cdot z + f(t, x,\mu,\alpha),
$$
and $\hat\alpha$ denotes its minimizer. 
Solutions of MFC can also be characterized through a similar FBSDE system, but in the latter case, the backward equation involves a partial derivative of the Hamiltonian with respect to the measure argument. See~\cite[Section 6.4.2]{MR3752669} for more details. Moreover, such FBSDE systems can also be obtained using dynamic programming, in which case $Y$ represents the value function instead of its gradient. 

\vskip 2pt
 All these systems are particular cases of the following general system of forward-backward SDEs of McKean-Vlasov type (MKV FBSDE for short), the system~\eqref{eq:MKV-FBSDE} being derived from the application of the Pontryagin stochastic maximum principle applied to our original MKV control problem
\begin{equation}
\label{eq:MKV-FBSDE-general}
\left\{
\begin{aligned}
	d X_t
	=
	& B\left(t, X_t, \cL(X_t), Y_t \right) dt
		+  \sigma d W_t,
	\\
	d Y_t
	=
	& - F\left(t, X_t, \cL(X_t), Y_t, \sigma^\dagger Z_t \right) dt
		+ Z_t d W_t,
\end{aligned}
\right.
\end{equation}
with initial condition $X_0$ having distribution $m_0$ and terminal condition $Y_T = G(X_T, \cL(X_T))$.

\vskip 2pt
The solution strategy is to replace the backward equation forced on us due to the optimization, by a forward equation and treat its initial condition, which is what we are looking for, as a control for a new optimization problem. This strategy has been successfully applied to problems in economic contract theory where it is known as Sannikov's trick. See for example \cite{MR1766421,MR2963805,MR3738664}.
Put it plainly, the strategy is a form of shooting method: the controller chooses the initial point and the volatility of the $Y$ process, and penalizes them proportionally to how far it is from matching the terminal condition. Specifically, we minimize over $y_0 : \RR^d \to \RR^d$ and $z: \RR_+ \times \RR^d \to \RR^{d \times d}$ the cost functional
\begin{equation}
\label{eq:cost-FBSDE-penalization}
	J_{FBSDE}(y_0, z) = \EE \left[ \, \left| Y^{y_0,z}_T - G(X^{y_0,z}_T, \cL(X^{y_0,z}_T)) \right|^2 \, \right]
\end{equation}
where $(X^{y_0,z}, Y^{y_0,z})$ solves
\begin{equation}
\label{eq:MKV-FBSDE-2forward}
\left\{
\begin{aligned}
d X^{y_0,z}_t
	=
	& B\left(t, X^{y_0,z}_t, \cL(X^{y_0,z}_t), Y^{y_0,z}_t \right) dt
		+  \sigma d W_t,
	\\
	d Y^{y_0,z}_t
	=
	& - F\left(t, X^{y_0,z}_t, \cL(X^{y_0,z}_t), Y^{y_0,z}_t, \sigma^\dagger z(t, X^{y_0,z}_t)  \right) dt
		+ z(t, X^{y_0,z}_t) d W_t,
\end{aligned}
\right.
\end{equation}
with \emph{initial} condition 
$X^{y_0,z}_0 = \xi_0 \in L^2(\Omega, \cF_0, \PP; \RR^d)$ and  
$Y^{y_0,z}_0 = y_0 (X_0)$. 
In some sense, the above problem is an optimal control problem of MKV dynamics if we view $(X^{y_0,z}_t,Y^{y_0,z}_t)$ as state and $(y_0,z)$ as control. It is rather special because the control is initial value and the volatility of the second component of the state, and looked for among feedback functions of the first component of the state. Under suitable conditions, the optimally controlled process $(X_t,Y_t)_t$ solves the FBSDE system~\eqref{eq:MKV-FBSDE-general} and vice-versa. 

In the same spirit as the method presented in Section~\ref{sec:directMethod}, we consider a finite-size population with $N$ particles and replace the controls $y_0$ and $z$ by neural networks, say $y_{0,\theta}$ and $z_{\omega}$ with parameters $\theta$ and $\omega$ respectively. We then discretize time with steps of size $\Delta t$. Let us denote by $\bJ^{N, \Delta t}_{FBSDE,S}(\theta,\omega) $ the analog of~\eqref{subchapRCML-num-eq:NN-MFC-cost-totalapprox-oneS} for the cost function~\eqref{eq:cost-FBSDE-penalization} stemming from the MKV FBSDE. Finally, we use SGD to perform the optimization. The method is summarized in Algorithm~\ref{subchapRCML-num-algo:SGD-FBSDE}. It is similar to Algorithm~\ref{subchapRCML-num-algo:SGD-MFC} so we only stress the main differences. The two neural networks could be taken with different architectures and their parameters optimized with different learning rates. In Section \ref{sec:numres} below, we illustrate the performance of this method on MKV FBSDEs coming from an MFG model of systemic risk.

\vskip 6pt
\begin{algorithm}[ht]
\DontPrintSemicolon
\KwData{An initial parameter $\theta^{(0)},\omega^{(0)}\in\Theta$; a number of steps $\mathtt{K}$;  a sequence $(\beta^{(\mathtt{k})})_{\mathtt{k}=0,\dots,\mathtt{K}-1}$ of learning rates.}
\KwResult{Parameters $(\theta,\omega)$ such that $(y_{0,\theta},z_{\omega})$ approximately minimizes $J_{FBSDE}$}
\Begin{
  \For{$\mathtt{k} = 0, 1, 2, \dots, \mathtt{K}-1 $}{
    Pick $S  = (\underline {\check X}_0, (\Delta \underline {\check W}_n)_{n})$\;
    Simulate $N$ trajectories for~\eqref{eq:MKV-FBSDE-2forward} with $y_0=y_{0,\theta^{(\mathtt{k})}}$ and $z = z_{\theta^{(\mathtt{k})}}$\;
    Compute the gradient $\nabla_{(\theta,\omega)} \bJ^{N, \Delta t}_{FBSDE,S}(\theta^{(\mathtt{k})},\omega^{(\mathtt{k})})$\;
    Set $(\theta^{(\mathtt{k}+1)}, \omega^{(\mathtt{k}+1)}) = (\theta^{(\mathtt{k})}, \omega^{(\mathtt{k})}) - \beta^{(\mathtt{k})} \nabla_{(\theta,\omega)} \bJ^{N, \Delta t}_{FBSDE,S}(\theta^{(\mathtt{k})})$ \;
    }
  \KwRet{ $(\theta^{(\mathtt{K})}, \omega^{(\mathtt{K})})$}
  }
\caption{SGD for MKV FBSDE\label{subchapRCML-num-algo:SGD-FBSDE}}
\end{algorithm}

\subsection{Numerical illustration: a toy model of systemic risk}
\label{sec:numres}
The following MFG model was introduced in~\cite{MR3325083} as an example which can be solved explicitly with a common noise, for finitely many players as well as in the mean field limit, in the open loop case as well as the closed loop set-up, and for which the master equation can be derived and solved explicitly. Individual players are financial institutions, and their states are the logarithms of their cash reserves. We assume that their evolutions are given by one dimensional diffusion processes involving a common noise $\bW^0$ and an independent idiosyncratic noise $\bW$. The costs take into account the rates of lending and borrowing and penalize departure from the aggregate state of the other institutions.

Because of the presence of the common noise the best response needs to be computed when the conditional flow of distributions is fixed. Due to the linear-quadratic nature of the model, the optimization is performed given the flow of the conditional mean log-monetary reserves $\bar {\mathbf m} = (\bar m_t)_{t \in [0,T]}$ which is adapted to the filtration generated by $\bW^0$. Assuming that the log-monetary reserve of a bank satisfies the SDE:
$$
	dX_t = [a (\bar m_t - X_t) + \alpha_t] dt + \sigma \left( \rho \, dW^0_t + \sqrt{1 - \rho^2} dW_t\right).
$$
where $a>0$ and $\rho \in [0, 1]$ builds dependence between the random shocks. The first term comes from the fact that the bank is assumed to borrow or lend to each other bank at a rate proportional to the difference between their log-monetary reserves. The term in $\alpha_t$ represents the rate at which the bank borrows or lends to a central bank. Each institution tries to minimize its expected cost given by:
\begin{align*}
	&J^{MFG}(\bar m, \alpha) 
	\\
	&= \EE\left[\int_0^T \left( \frac{1}{2}\alpha_t^2 - q \alpha_t (\bar m_t - X_t) + \frac{\epsilon}{2} (\bar m_t - X_t)^2 \right) dt + \frac{c}{2} (\bar m_T - X_T)^2\right]
\end{align*}
$q$, $\epsilon$, $c$ and $\sigma$ being positive constants satisfying $q \le \epsilon^2$ so that the running cost is jointly convex in the state and the control variables.  Here, $q$ can be interpreted as a parameter chosen by a regulator to incentivize borrowing or lending: if the log-monetary reserve $X_t$ of the bank is smaller than the average $\bar m_t$, then the bank has an incentive to choose a positive control $\alpha_t$, meaning that it borrows; similarly, if $X_t > \bar m_t$, then the banks has an incentive to choose $\alpha_t>0$. 
For more details on the interpretation of this model in terms of systemic risk, the reader is referred to~\cite{MR3325083}. 
 The model is of linear-quadratic type and hence has an explicit solution through a Riccati equation, which we use as a benchmark for comparison with our numerical results. 

\vskip 2pt

Note that since this example is a MFG, we can not use Algorithm~\ref{subchapRCML-num-algo:SGD-MFC} to compute directly the equilibrium. Instead, we use Algorithm~\ref{subchapRCML-num-algo:SGD-FBSDE} to solve the appropriate FBSDE system (we omit this system here for brevity; see~\cite{MR3325083}). In order to deal with the additional randomness induced by the common noise, we add $\bar m_t$ as an argument of the neural networks playing the roles $y_0(\cdot)$ and $z(\cdot)$ introduced above.

\begin{figure}[ht]
\centering
\begin{subfigure}{.45\textwidth}
  \centering
  \includegraphics[width=\linewidth]{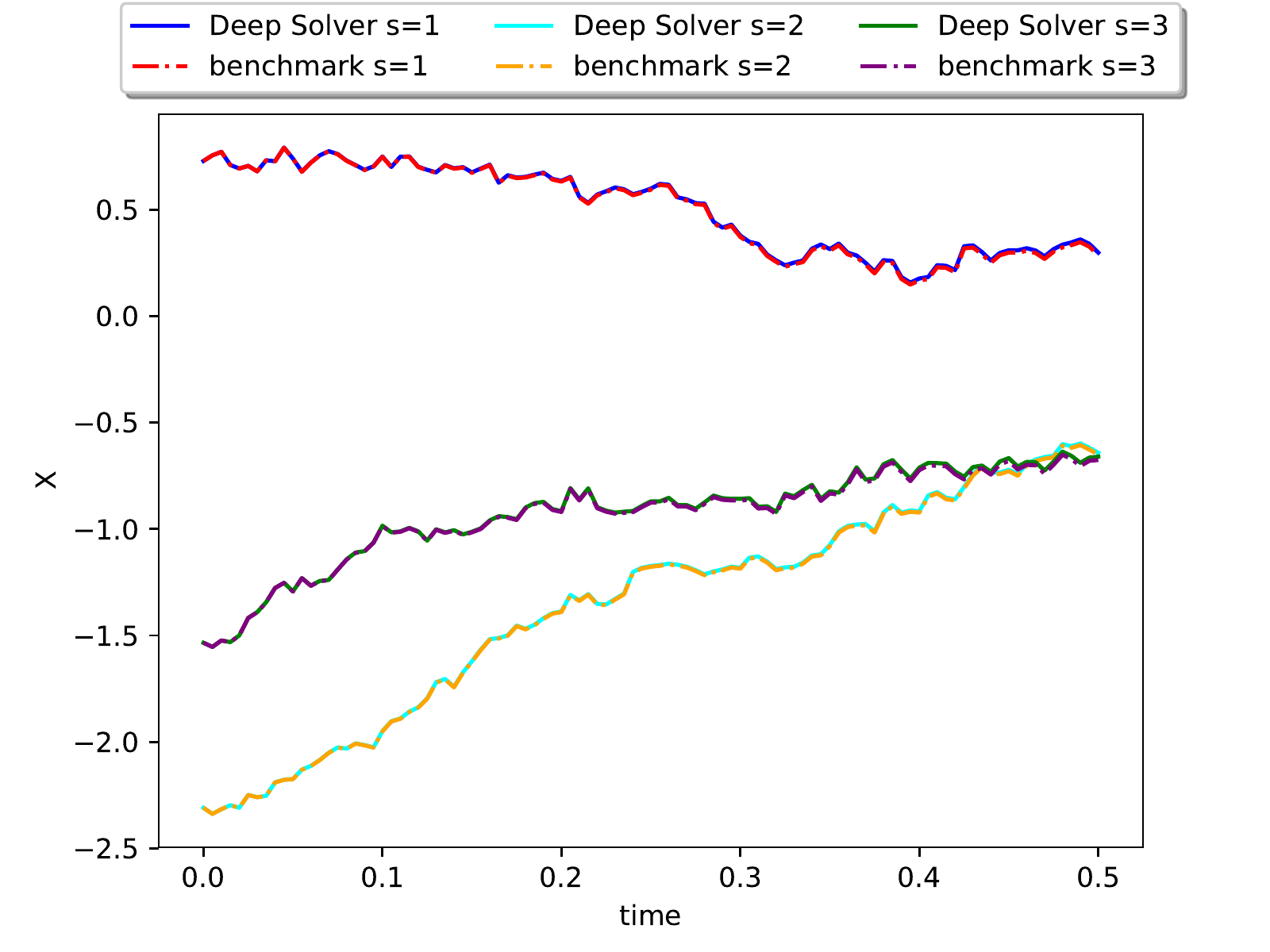}
  \caption*{Trajectory of $X^i, i=1,2,3$}
\end{subfigure}
\begin{subfigure}{.45\textwidth}
  \centering
  \includegraphics[width=\linewidth]{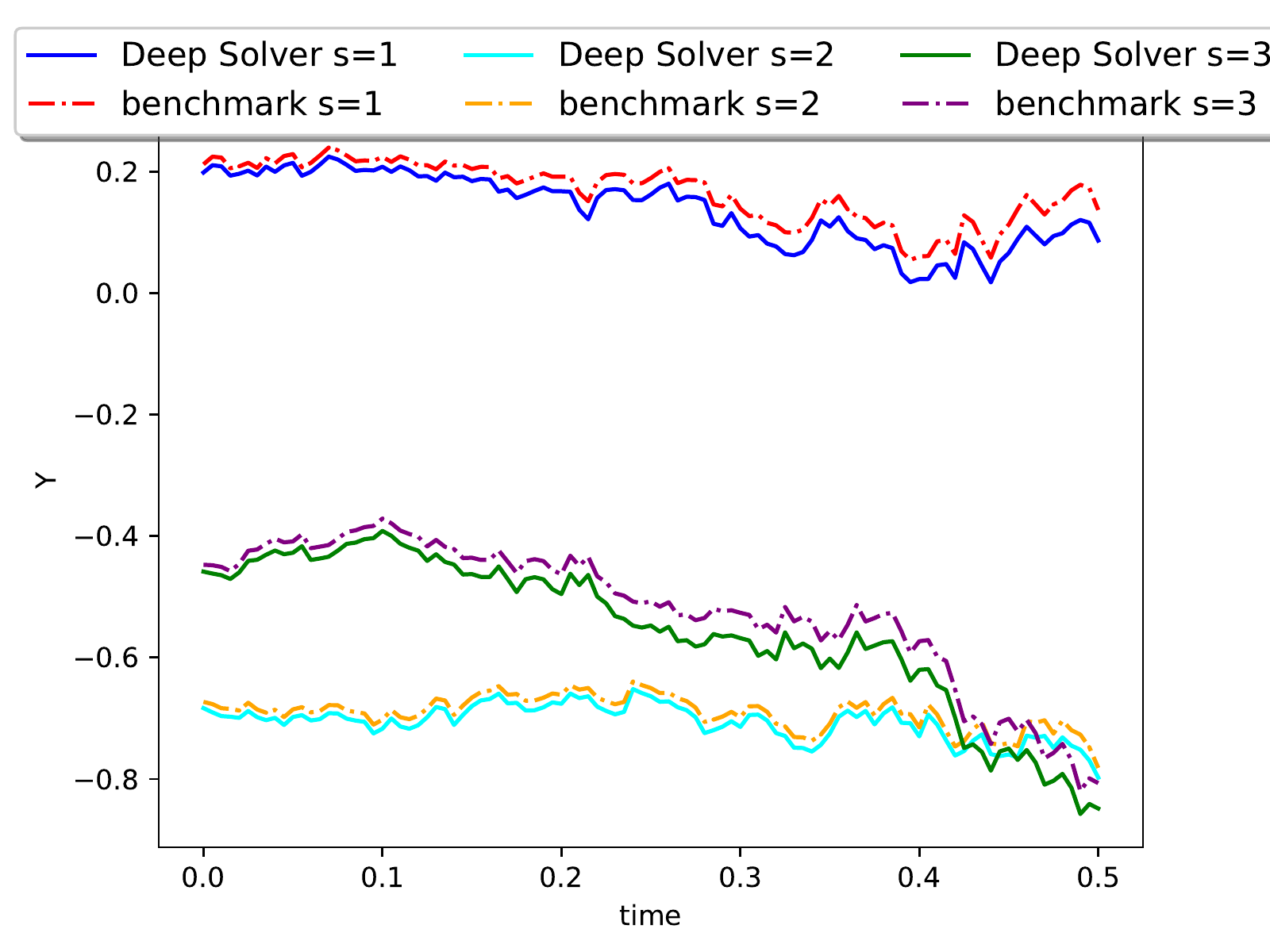}
  \caption*{Trajectory of $Y^i, i=1,2,3$}
\end{subfigure}
\caption{Systemic risk MFG example solved by Algorithm~\ref{subchapRCML-num-algo:SGD-FBSDE}. Three sample trajectories: solution computed by deep solver (full lines, in cyan, blue, and green) and by analytical formula (dashed lines, in orange, red and purple).}
\label{fig:ex-mfg-sysrisk-traj}
\end{figure}

Figure~\ref{fig:ex-mfg-sysrisk-traj} displays sample trajectories of $X^i$ and $Y^i$ for three different values of $i$. We can see that the approximation is better for $X^i$ than for $Y^i$, particularly towards the end of the time interval. This is probably due to the fact that the latter is supposed to solved a BSDE but we replaced it by a forward equation so errors accumulate along time. However the error decreases with the number of time steps, particles and units in the neural network. See~\cite{carmona2019convergence2} for more details. For the numerical tests presented here, we used $\sigma = 0.5, \rho = 0.5, q = 0.5, \epsilon = q^2 + 0.5 = 0.75, a = 1, c = 1.0$ and $T = 0.5$.

\section{DGM method for mean field PDEs}
\label{sec:MFDGM}

In this section we present an adaptation of the Deep Galerkin Method (DGM) introduced in~\cite{MR3874585} for MFGs and MFC. In this context, it can be used to solve the forward-backward PDE system~\cite{al2018solving,carmona2021convergence1,ruthotto2020machine,cao2020connecting,lin2020apac,AMSnotesLauriere} or some forms of the Master equation~\cite{AMSnotesLauriere}. The key idea is to replace the unknown function by a neural network and to tune the parameters to minimize a loss function based on the residual of the PDE. After presenting the main ideas, we illustrate this method on a model of optimal execution.

\subsection{Description of the method}

We present the method with the example of a finite horizon MFG model on a compact domain. We rewrite the MFG PDE system as a minimization problem where the control is the pair of density and value function, and the loss function is the sum of the PDE residuals and terms taking into account the boundary conditions. The same approach can be adapted to the ergodic setting, where initial and terminal conditions are replaced by normalization conditions see~\cite{carmona2021convergence1}. 
In a finite horizon MFG as defined in Definition~\ref{def:MFGeq}, the optimal control can be characterized (under suitable conditions, see~\cite{MR2295621}), as:
\begin{equation*}
	\hat{\ctrl}(x, m(t), \grad u(t,x)) = 
	{\rm{argmin}}_{a \in \RR^k}  \Bigl( f(x, m(t), a)+ \grad u(t,x) \cdot b(x, m(t), a)\Bigr),
\end{equation*}

where $(m, u)$ solve the
forward-backward PDE system on $\domT$:
\begin{empheq}[left=\empheqlbrace]{align}
\label{eq:MFG-KFP}
	&0 =
	\displaystyle \partial_t m (t,x)  - \nu \Delta m(t,x) + \div\Bigl( m(t,x) \partial_q H^*( x, m(t),  \grad u(t,x))\Bigr)
	\\
\label{eq:MFG-HJB}
&0 =
	\partial_t u (t,x) + \nu \Delta u(t,x) + H^*( x, m(t), \grad u(t,x))
\end{empheq}
with the initial and terminal conditions:
\begin{align*}
  	& m(0,x) = m_0(x), \quad u(T,x) = g(x, m(T)), \qquad x \in \dom,
\end{align*}
where $\nu = \tfrac{\sigma^2}{2}$ and $H^*$ is the optimized Hamiltonian, defined as:
\begin{align*}
	H^*( x, m, q) := {\rm{min}}_{a \in \RR^k}  \Bigl( f(x, m, a)+ q \cdot b(x, m, a)\Bigr).
\end{align*}
The Kolmogorov-Fokker-Planck (KFP) equation~\eqref{eq:MFG-KFP} describes the evolution of the population distribution, while the Hamilton-Jacobi-Bellman (HJB) equation~\eqref{eq:MFG-HJB} describes the evolution of the value function. It is obtained for instance by a dynamic programming argument for the infinitesimal player's optimization problem given the flow of population density. These PDEs are coupled, hence we can not solve one before the other one.

\vskip 2pt

For simplicity, we replace the domain $\dom$ by a compact subset $\tilde\dom$. We denote $\tilde\domT = [0,T] \times \tilde\dom$. We introduce the following loss function:
\begin{equation}
\label{eq:loss-DGM-total}
	L(m,u) = L^{(KFP)}(m,u) + L^{(HJB)}(m,u)
\end{equation}
where
\begin{equation}
\label{eq:loss-DGM-KFP}
\begin{split}
	 L^{(KFP)}(m,u)
	 &= C^{(KFP)} \left\| \partial_t m  - \nu \Delta m + \div\Bigl( m \partial_q H^*( x, m(t),  \grad u)\Bigr) \right\|_{L^2(\tilde\dom)}
	 \\
	 & \qquad + C^{(KFP)}_0 \left\| m(0) - m_0 \right\|_{L^2(\tilde\dom)}
\end{split}
\end{equation}
and 
\begin{equation}
\label{eq:loss-DGM-HJB}
\begin{split}
	 L^{(HJB)}(m,u)
	&= C^{(HJB)} \left\| \partial_t u  + \nu \Delta u + H^*( \cdot, m(t), \grad u) \right\|_{L^2(\tilde\domT)}
	 \\
	 & \qquad + C^{(HJB)}_T \left\| u(T) - g(\cdot, m(T)) \right\|_{L^2(\tilde\dom)}.
\end{split}
\end{equation}
Each component of the loss $L$ in~\eqref{eq:loss-DGM-total} encodes one of the two PDEs of the optimality system~\eqref{eq:MFG-KFP}--\eqref{eq:MFG-HJB} with  one term for the PDE residual and one term for the initial or terminal condition. The positive constants $C^{(KFP)}, C^{(KFP)}_0, C^{(HJB)},$ and $C^{(HJB)}_T$ give more or less importance to each component.  On a bounded domain with boundary condition, more penalty terms could be included. Note that $L(m,u) = 0$ if $(m,u)$ is a smooth enough solution to the PDE system~\eqref{eq:MFG-KFP}--\eqref{eq:MFG-HJB}. 
A similar system can be derived to characterize the optimal control of a MFC and the method can be adapted to this setting. See~\cite{carmona2021convergence1}.

\vskip 6pt
Replacing $m$ and $u$ by neural networks is the lynchpin of the algorithm. We denote by $m_{\theta_1}$ and $u_{\theta_2}$ these neural nets parameterized by $\theta_1$ and $\theta_2$ respectively. As in the method discussed in the previous sections, the integrals on $\tilde\domT$ (resp. $\tilde\dom$) are interpreted as expectations with respect to a uniform random variable over $\tilde\domT$ (resp. $\tilde\dom$), and we use SGD to minimize the total loss function $L$. More precisely, for a given $\mathbf{S} = (S, S_0, S_T)$ where $S $ is a  finite set of points in $\tilde\domT$, and $S_0 $ and $S_T$ are finite sets of points in $\dom$, we define the empirical loss function as 
\begin{equation}
\label{eq:loss-DGM-system}
	L_{\mathbf{S}}(\theta) = 
	L_{\mathbf{S}}^{(KFP)}(\theta) + L_{\mathbf{S}}^{(HJB)}(\theta),
\qquad\qquad\theta = (\theta_1,\theta_2)
\end{equation}
where
\begin{align*}
	L_{\mathbf{S}}^{(KFP)}(\theta) 
	= \, & C^{(KFP)}\Big(\frac{1}{|S|}\sum_{(t,x) \in S} \Big|\partial_t m_{\theta_1}(t,x)  - \nu \Delta m_{\theta_1}(t,x) 
	\\
	&\qquad\qquad\qquad + \div\Bigl( m_{\theta_1}(t,x) \partial_q H^*( x, m_{\theta_1}(t),  \grad u_{\theta_2}(t,x))\Bigr) \Big |^2\Big)^{1/2}
	\\
	&\qquad +  C^{(KFP)}_0 \left(\frac{1}{|S_0|}\sum_{x \in S_0} \left | m(0,x) - m_0(x) \right |^2\right)^{1/2}
\end{align*}
and
\begin{align*}
	L_{\mathbf{S}}^{(HJB)}(\theta) 
	= \, & C^{(HJB)} \left(\frac{1}{|S|}\sum_{x \in S} \left | \partial_t u_{\theta_2}(t,x)  + \nu \Delta u_{\theta_2}(t,x) + H^*( x, m_{\theta_1}(t), \grad u_{\theta_2}(t,x))\right |^2 \right)^{1/2}
	\\
	&\qquad + C^{(HJB)}_T \left(\frac{1}{|S_T|}\sum_{x \in S_T} \left | u_{\theta_2}(T,x) - g(x, m_{\theta_1}(T)) \right |^2\right)^{1/2}.
\end{align*}

The method is summarized in Algorithm~\ref{subchapRCML-num-algo:DGM-MFG}. The two neural networks could be taken with different architectures and their parameters optimized with different learning rates. The convergence of the neural network approximation was discussed in~\cite{MR3874585} in the context of a single PDE using a standard universal approximation theorem. Unfortunately, this does not shed any light on the rate of convergence. A rate of convergence can be obtained by using more constructive approximation results with neural networks. See~\cite{carmona2021convergence1} and the references therein. In turn, this property leads to bounds on both the loss function of the algorithm and the error on the value function of the control problem. However, to the best of our knowledge, the convergence of the algorithm towards approximately optimal parameters remains to be proved, as for the other methods presented in this subchapter.

\vskip 6pt
\begin{algorithm}[ht]
\DontPrintSemicolon
\KwData{Initial parameters $\theta^{(0)} = (\theta^{(0)}_1,\theta^{(0)}_2)\in\Theta$; a number of steps $\mathtt{K}$;  a sequence $(\beta^{(\mathtt{k})})_{\mathtt{k}=0,\dots,\mathtt{K}-1}$ of learning rates.}
\KwResult{Parameters $\theta = (\theta_1,\theta_2)$ such that $(m_{\theta_1}, u_{\theta_2})$ approximately minimize $L$ defined in~\eqref{eq:loss-DGM-total}}
\Begin{
  \For{$\mathtt{k} = 0, 1, 2, \dots, \mathtt{K}-1 $}{
    Pick $\mathbf{S} = (S, S_0, S_T)$\;
    Compute the gradient $\nabla_\theta L_{\mathbf{S}}(\theta^{(\mathtt{k})})$, see~\eqref{eq:loss-DGM-system}\;
    Set $\theta^{(\mathtt{k}+1)} = \theta^{(\mathtt{k})} - \beta^{(\mathtt{k})} \nabla_\theta L_{\mathbf{S}}(\theta^{(\mathtt{k})})$ \;
    }
  \KwRet{ $\theta^{(\mathtt{K})}$}
  }
\caption{DGM for MFG PDE system\label{subchapRCML-num-algo:DGM-MFG}}
\end{algorithm}

\vskip 6pt
An important advantage of the DGM method is its flexibility and its generality: in principle, it can be applied to almost any PDE since it is agnostic to the structure of the PDE in question, or in the extension described above, of the PDE system. In tailoring the strategy to the specifics of our system, our main challenge was the choice of the relative weights to be assigned to the various terms in the aggregate loss function. If they are not chosen appropriately, SGD can easily be stuck in local minima. For example, if the weights $C^{(KFP)}_0$ and $C^{(HJB)}_T$ are not large enough, the neural networks might find trivial solutions minimizing the residuals while ignoring the initial and terminal conditions. However, if these weights are too large, the neural networks might satisfy very well these conditions without solving very precisely each PDE on the interior of the domain. See~\cite{carmona2021convergence1} for a more detailed discussion on this aspect.

\subsection{Numerical illustration: a crowded trade model}
This is a model of optimal execution similar to the one studied in Subsection \ref{sub:price_impact}.
Here, we follow~\cite{MR3805247} and we assume that a broker is instructed by a client to liquidate $Q_0$ shares of a specific stock by a fixed time horizon $T$, and that this broker is representative of a large population of brokers trying to do exactly the same thing. Such a situation occurs when a large number of index trackers need suddenly to rebalance their portfolios because the composition of the market index they track is changed. So our typical broker tries
to maximize the quantity:
$$
    \mathbb{E} \left[ X_T + Q_T(S_T - A Q_T) - \phi \int_0^T |Q_t|^2 dt \right]
$$
where at time $t\in[0,T]$, $S_t$ is the price of the stock, $Q_t$ is the inventory (i.e. number of shares) held by the broker,  and $X_t$ is their wealth. The constant $\phi>0$ weigh a penalty for holding inventory through time while $A>0$ plays a similar role at the terminal time. The dynamics of these three state variables are given by:  
subject to:
\begin{equation*}
\begin{cases}
    dS_t = \gamma \bar{\mu}_t  dt + \sigma dW_t\\
    dQ_t = \ctrl_t dt\\
    dX_t = -\ctrl_t(S_t+\kappa \ctrl_t)dt.
\end{cases}
\end{equation*}
The time evolution of the price $S_t$ is subject to random shocks with standard deviation $\sigma$ where the innovation $dW_t$ is given by the increments of a standard Brownian motion, and a drift accounting for a permanent price impact $\gamma\bar{\mu}_t$ resulting from the aggregate trading rate $\bar{\mu}_t$ of all the brokers multiplied by a constant $\gamma>0$.
The rate of trading $\ctrl_t$ is the control of the broker.
Finally, the constant $\kappa>0$ account for a quadratic transaction cost.
 Except for the fact that $\bar{\mu}_t$ is here endogenous, this is the model considered in \cite{MR3500455}, to which a deep learning method has been applied in~\cite{leal2020learning} to approximate the optimal control on real data.

\begin{remark}
The current model has two major differences with the model considered earlier in Subsection \ref{sub:price_impact}. It does not belong to the class of linear-quadratic models because the transaction costs entering the dynamics are quadratic in the control. But most importantly, the broker's inventory does not have a Brownian component. The presence of a quadratic variations term in the dynamics of the inventory was demonstrated in \cite{CarmonaWebster} running econometric tests on high frequency market data. This was one of the reasons for the choice of the model used in Subsection \ref{sub:price_impact}. Surprisingly, it is shown in \cite{CarmonaLeal} that the inclusion of a Brownian motion component in the dynamics of the inventory process $Q_t$ does not require significant changes to the proof, including the form of the ansatz for the value function.
\end{remark}

In any case, when the flow $(\bar{\mu}_t)_{0\le t\le T}$ is fixed, the optimization problem involved in the computation of the best response reduces to an HJB equation whose solution $V(t,x,s,q)$ can be found like in \cite{MR3500455} by formulating the ansatz 
 $V(t,x,s,q) = x + qs + v(t,q)$ for some function $v$. Rewriting the HJB equation one sees that $v$ must solve the equation:
$$
    -\gamma\bar{\mu} q = \partial_t v - \phi q^2 + \sup_\ctrl \{\ctrl \partial_q v - \kappa \ctrl^2\}
$$
with terminal condition $v(T,q) = - A q^2$, the optimal control being $\ctrl^*_t(q) = \frac{\partial_q v(t,q)}{2\kappa}$. Accordingly, if we denote by $m(t,\cdot)$ the distribution of inventories at time $t$, the aggregate trading rate is given by:
$$
    \bar{\mu}_t = \int \ctrl^*_t(q) m(t,dq) = \int \frac{\partial_q v(t,q)}{2\kappa} m(t, dq),
$$
in equilibrium since we use the optimal control. Since the evolution of the inventory distribution can be captured by the Kolmogorov-Fokker-Planck partial differential equation:
$$
    \partial_t m + \partial_q\left(m \frac{\partial_q v(t,q)}{2\kappa}\right) = 0,
$$
with a given initial condition $m(0,\cdot) = m_0$, the solution of the MFG can be characterized by the PDE system is (see \cite{MR3805247} for more details):
\begin{equation}
\label{eq:CL-PDE-reduced}
\left\{
\begin{aligned}
    &\quad -\gamma\bar{\mu} q = \partial_t v - \phi q^2 + \frac{|\partial_q v(t,q)|^2}{4\kappa}
    \\
    &\quad \partial_t m + \partial_q\left(m \frac{\partial_q v(t,q)}{2\kappa}\right) = 0
    \\
    &\quad \bar{\mu}_t = \int \frac{\partial_q v(t,q)}{2\kappa} m(t, dq)
    \\
    &\quad m(0,\cdot) = m_0, v(T,q) = - A q^2. 
\end{aligned}
\right.
\end{equation}
Note that the mean field interactions are through $\bar{\mu}_t$, which is a non-local (in space) term involving the derivative of the solution to the HJB equation.  

This PDE system has been solved with the DGM method in~\cite{al2018solving} after a change of variable for the distribution. Here, for the sake of numerical illustration, we present results based on directly solving this system by following the methodology discussed above, suitably modified for the time-dependent PDE system~\eqref{eq:CL-PDE-reduced}. The initial and terminal conditions are imposed by penalization. The non-local term is estimated with Monte Carlo samples. For the results presented here, we used the following values for the parameters: $T=1$, $\sigma=0.3$, $A = 1$, $\phi = 1$, $\kappa = 1$, $\gamma = 1$, and a Gaussian initial distribution with mean $4$ and variance $0.3$.

\begin{figure}[ht]
\centering
\begin{subfigure}{.45\textwidth}
  \centering
  \includegraphics[width=\linewidth]{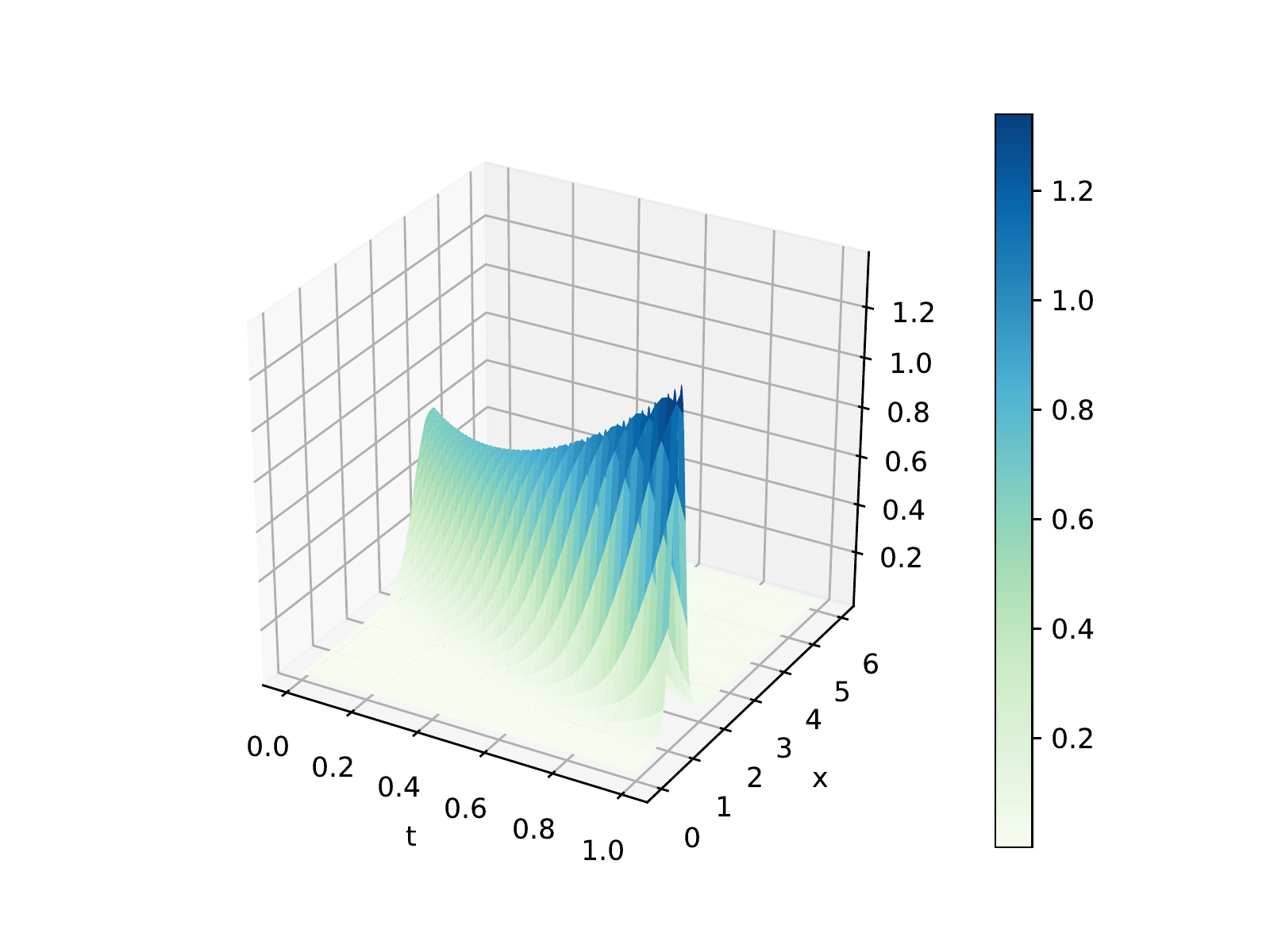}
\end{subfigure}
\begin{subfigure}{.45\textwidth}
  \centering
  \includegraphics[width=\linewidth]{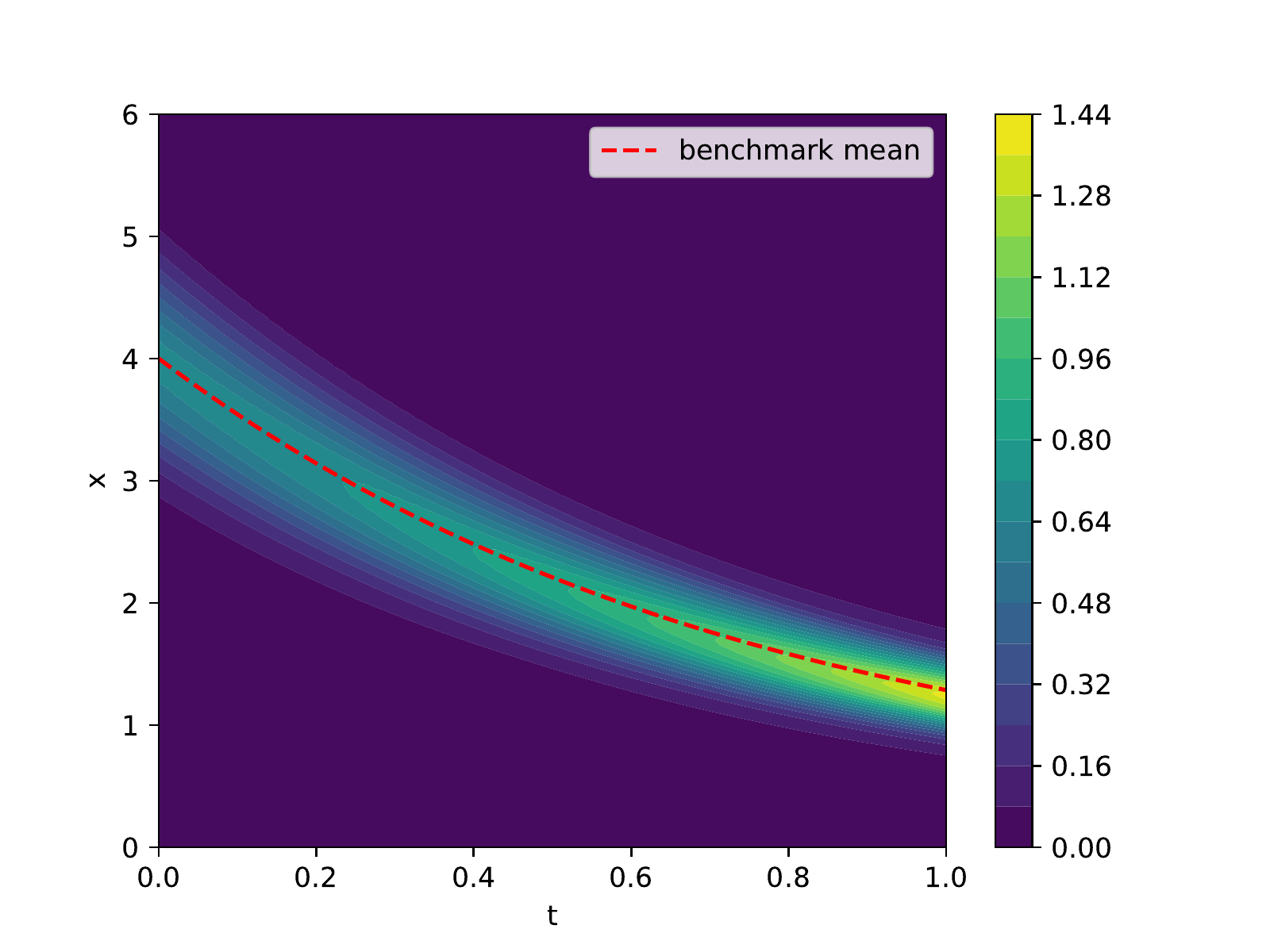}
\end{subfigure}
\caption{Trade crowding MFG example solved by Algorithm~\ref{subchapRCML-num-algo:DGM-MFG}. Evolution of the distribution $m$: surface (left) and contour (right).}
\label{fig:ex-mfg-crowd-trade-distrib}
\end{figure}

\begin{figure}[ht]
\centering
\begin{subfigure}{.3\textwidth}
  \centering
  \includegraphics[width=\linewidth]{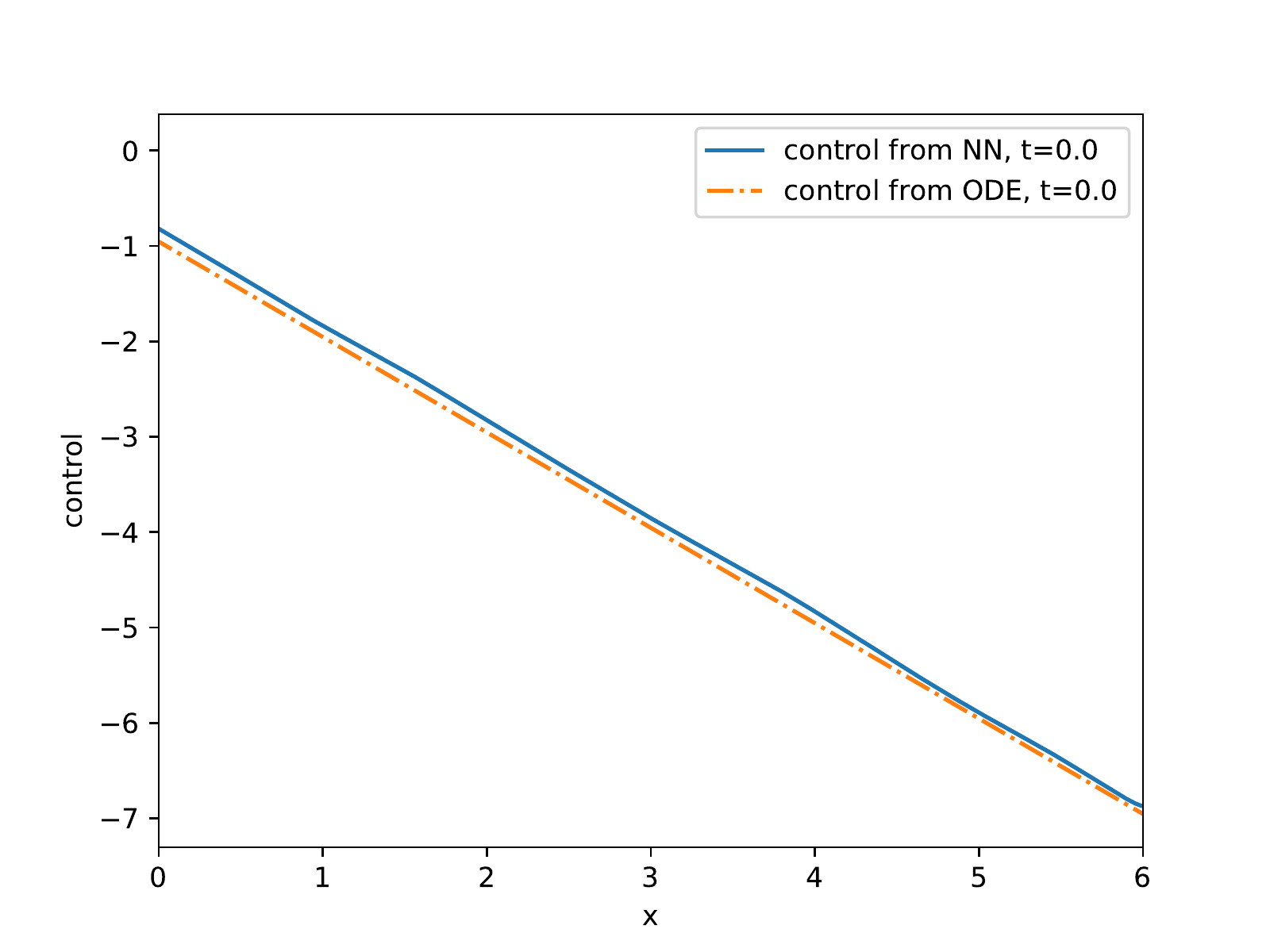}
\end{subfigure}
\begin{subfigure}{.3\textwidth}
  \centering
  \includegraphics[width=\linewidth]{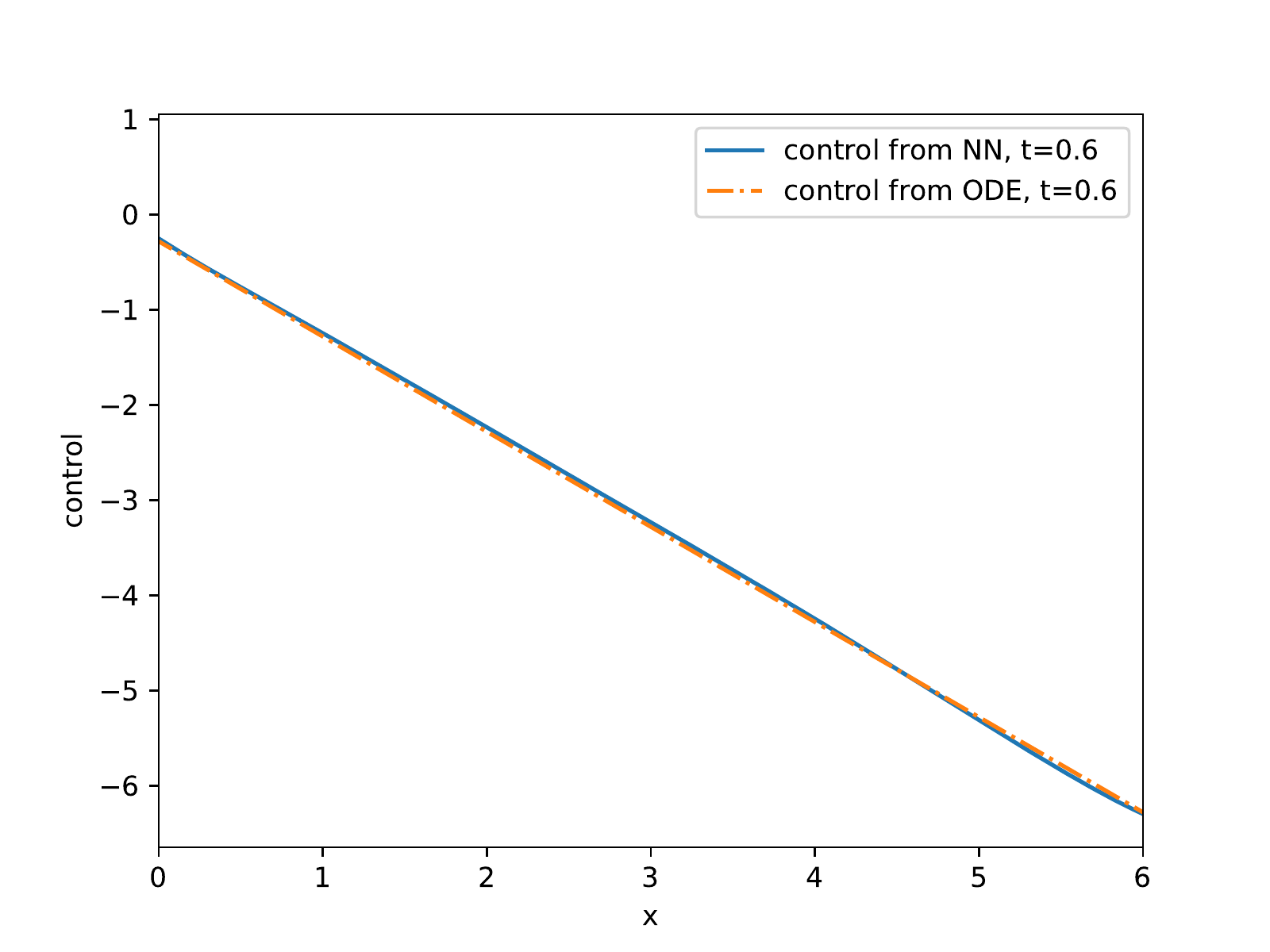}
\end{subfigure}
\begin{subfigure}{.3\textwidth}
  \centering
  \includegraphics[width=\linewidth]{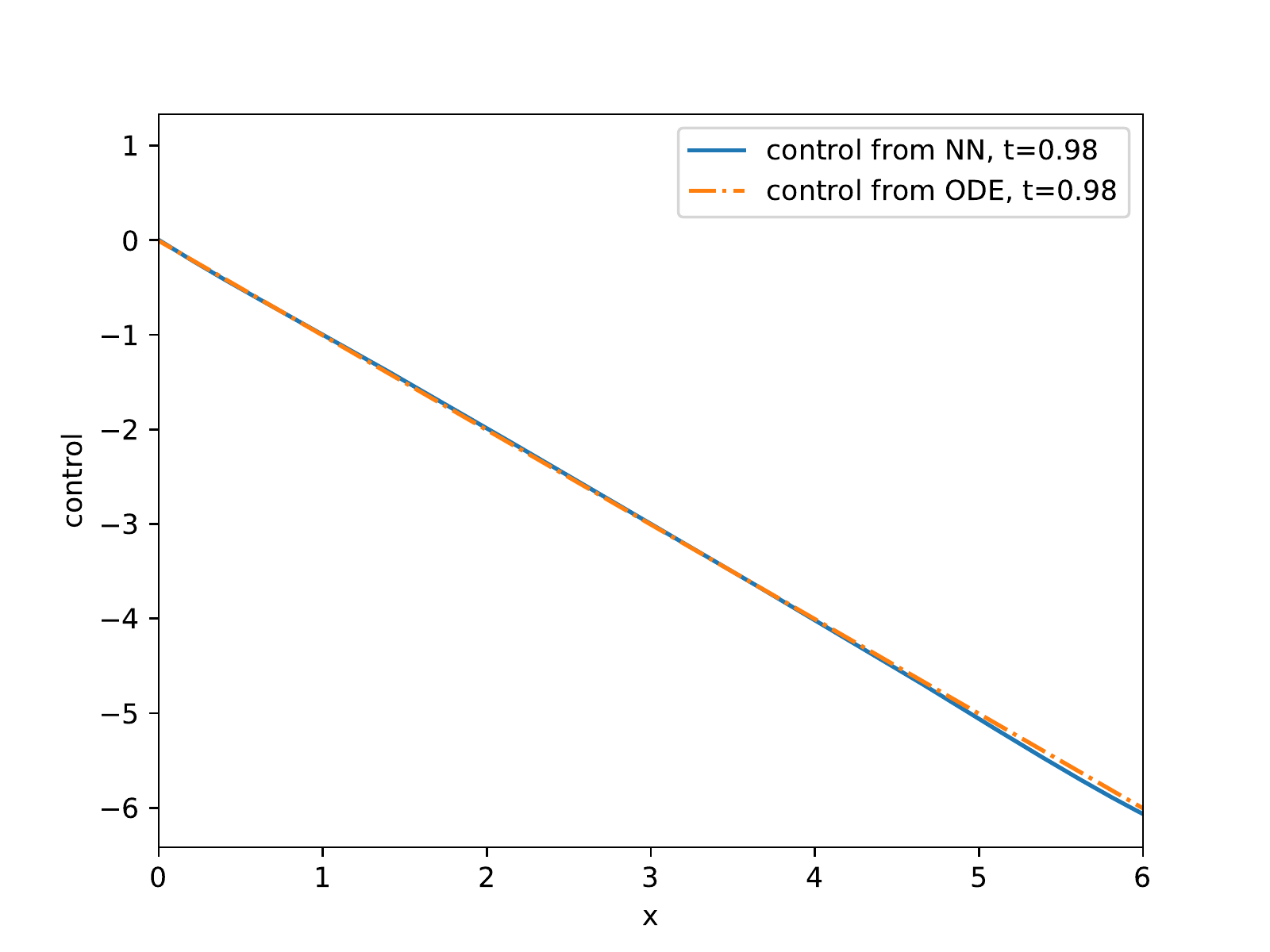}
\end{subfigure}
\caption{Trade crowding MFG example solved by Algorithm~\ref{subchapRCML-num-algo:DGM-MFG}. Optimal control $\ctrl^*$ (dashed line) and learnt control (full line)  at three different time steps.}
\label{fig:ex-mfg-crowd-trade-value}
\end{figure}

The evolution of the distribution $m$ is displayed in Figure~\ref{fig:ex-mfg-crowd-trade-distrib} while the value function $v$ and the optimal control $\ctrl^*$ are displayed in Figure~\ref{fig:ex-mfg-crowd-trade-value}. As expected from the theory, the distribution concentrates close to $0$ and we recover a linear control, which matches the optimal one obtained with semi-explicit formula (see~\cite{MR3805247} for more details).  For the neural network approximating the density, on the last layer, we used an exponential activation function. This ensures that the density is always non-negative.

\section{Conclusion}
\label{sec:ccl}

In this subchapter we have presented three families of strategies to solve MFC and MFG. The first one is designed to solve MFC problems by directly trying to minimize the cost functional after replacing the control function by a neural network. The second approach tackles generic MKV FBSDE systems and uses a (stochastic) shooting method where the unknown starting point $Y_0$ of the backward variable and the $Z$ component are learnt as neural network functions of the state. The last approach solves mean field PDE systems by minimizing the residuals when the unknown functions are replaced by neural networks. We have illustrated these methods on stylized models arising in finance. The expressive power of neural networks let us expect that these methods will allow researchers and practitioners to solve much more complex and realistic models. The development of sample efficient methods able to learn solutions on real data while taking into account mean field interactions seem particularly relevant for future applications.

\bibliographystyle{abbrv}
\bibliography{chaptermlmfg}

\end{document}